\documentclass[12pt,reqno]{amsart}

\usepackage{amssymb,amsmath, amsfonts, latexsym}
\newtheorem{thm}[equation]{Theorem}
\newtheorem{lemma}[equation]{Lemma}
\newtheorem{prop}[equation]{Proposition}
\newtheorem{cor}[equation]{Corollary}

\theoremstyle{definition}
\newtheorem{definition}[equation]{Definition}
\newtheorem{recipe}[equation]{Recipe}

\theoremstyle{remark}
\newtheorem{remark}[equation]{Remark}

\numberwithin{equation}{section}

\begin{document}
\newcommand{\R}{\mathbb{R}}
\renewcommand{\P}{\mathbb{P}}
\newcommand{\E}{\mathbb{E}}
\newcommand{\1}{\textbf{1}}

\newcommand{\eqdef}{\stackrel{\vartriangle}{=}}
\newcommand{\dom}{\mathrm{dom\ }}
\newcommand{\icordom}{\mathrm{icordom\ }}
\newcommand{\cl}{\mathrm{cl\ }}
\newcommand{\inter}{\mathrm{int\ }}

\newcommand{\logsup}{\limsup_{n\rightarrow \infty}\frac 1n\log}
\newcommand{\loginf}{\liminf_{n\rightarrow \infty}\frac 1n\log}

\newcommand{\lsc}{lower semicontinuous}

\newcommand{\X}{\mathcal{X}}
\newcommand{\Y}{\mathcal{Y}}
\newcommand{\IX}{\int_\mathcal{X}}
\newcommand{\N}{\mathcal{N}}
\newcommand{\PX}{\mathcal{P}(\X)}
\newcommand{\NF}{\mathcal{P}_{\mathcal{F}}}


\newcommand{\pf}{(\psi, \varphi)}
\newcommand{\mc}{\circledast}
\newcommand{\Hnm}{H(\nu\mid\mu)}


\title{A large deviation approach to some transportation cost inequalities}

\author{Nathael Gozlan}
\address{(Nathael Gozlan) Modal-X, Universit\'e Paris 10. B\^at. G, 200 av. de la R\'epublique. 92001 Nanterre Cedex, France}
\email{nathael.gozlan@u-paris10.fr}

\author{Christian L\'eonard}
\address{(Christian L\'eonard) Modal-X, Universit\'e Paris 10. B\^at. G, 200 av. de la R\'epublique. 92001 Nanterre Cedex, France}
\address{(Christian L\'eonard) CMAP, \'Ecole Polytechnique. 91128 Palaiseau Cedex, France}
\email{christian.leonard@polytechnique.fr}

\date{July 2005}



\begin{abstract}
New transportation cost inequalities are derived by means of
elementary large deviation reasonings. Their dual characterization
is proved; this provides an extension of a well-known result of
S.~Bobkov and F.~G\"otze. Their tensorization properties are
investigated. Sufficient conditions (and necessary conditions too)
for these inequalities are stated in terms of the integrability of
the reference measure. Applying these results leads to new
deviation results: concentration of measure and deviations of
empirical processes.
\end{abstract}

\maketitle


\section{Introduction}

In the whole paper, $\X$ is a Polish space equipped with its Borel
$\sigma$-field. We denote $\PX$ the set of all probability
measures on $\X.$

\subsection{Transportation cost inequalities and concentration of measure}
Let us first recall what  transportation cost inequalites are and
their well known consequences in terms of concentration of
measure.

\par\medskip\noindent
 \textbf{Transportation cost.}\  Let $c: \X\times\X\rightarrow
[0,\infty)$ be a measurable function on the product space
$\X\times\X.$ For any couple of probability measures $\mu$ and
$\nu$ on $\X,$ the transportation cost (associated with the cost
function $c$) of $\mu$ on $\nu$ is
\begin{equation*}
\mathcal{T}_c(\mu, \nu)=\inf_\pi \int_{\X\times\X}
c(x,y)\,\pi(dxdy)\in [0,\infty]
\end{equation*}
where the inf is taken over all probability measures $\pi$ on
$\X\times\X$ with first marginal $\pi(dx\times\X)=\mu(dx)$ and
second marginal $\pi(\X\times dy)=\nu(dy).$

\par\medskip\noindent
 \textbf{$T_p$-inequalities.}\
Popular cost functions are $c(x,y)=d(x,y)^p$ where $d$ is a metric
on $\X$ and $p\geq 1.$ It is known that for some $\mu\in\PX$ and
$p\geq 1$ one can prove the following \emph{transportation cost
inequality}
\begin{equation}\label{eq-Tp}
\mathcal{T}_{d^p}(\mu,\nu)^{1/p} \leq \sqrt{2C
H(\nu\mid\mu)},\quad \forall\nu\in\PX
\end{equation}
for some positive constant $C,$ where
  $H(\nu\mid\mu)$ is the
\emph{relative entropy} of $\nu$ with respect to $\mu$ defined by
\[
H(\nu\mid\mu)=\int_\X \log\left(\frac{d\nu}{d\mu}\right)\,d\nu
\]
if $\nu$ is absolutely continuous with respect to $\mu$ and
$H(\nu\mid\mu)=\infty$ otherwise. In presence of the family of
inequalities (\ref{eq-Tp}), one says that $\mu$ satifies $T_p(C).$
\\
For instance, Csisz\'ar-Kullback-Pinsker's inequality, see
(\ref{eq-Pinsker}), is $T_1(1)$ with the Hamming's metric
$d(x,y)=\1_{x\not = y.}$ Csisz\'ar-Kullback-Pinsker's inequality
is often called Pinsker's inequality, it will be refered later as
CKP inequality. It holds for any $\mu\in\PX.$ On the other hand,
$T_2$-inequalities are much more difficult to obtain. It is shown
in the articles by F.~Otto and C.~Villani \cite{OVill00} and by
S.~Bobkov, I.~Gentil and M.~Ledoux \cite{BGL01}, that if $\mu$
satisfies the logarithmic Sobolev inequality, then it also
satisfies $T_2.$ A standard example of probability measure $\mu$
that satisfies $T_2$ is the normal law. In \cite{Tal96a},
M.~Talagrand has given a proof of $T_2(C)$ for the standard normal
law not relying on any log-Sobolev inequality, for the sharp
constant $C=1.$

\par\medskip\noindent
 \textbf{Concentration of measure.} \
As a consequence of $T_1(C),$ K.~Marton \cite{Mar86, Mar96} has
obtained the following \emph{concentration inequality} for $\mu:$
\begin{equation}\label{eq-10}
\mu(\{x; d(x,A)>r\})\leq
\exp\left[-\left(\frac{r}{\sqrt{2C}}-\sqrt{\log 2}\right)^2\right]
\end{equation}
for all measurable subset $A$ such that $\mu(A)\geq 1/2$ and all
$r\geq \sqrt{2C\log 2}.$ \emph{Marton's concentration argument}
easily extends to more general situations. This is of considerable
importance and justifies the search for $T_1$-inequalities.

\par\medskip\noindent
 \textbf{Product of measures.}\
Suppose that $\mu_1,\dots,\mu_n$ satisfy respectively
$T_p(C_1),\dots, T_p(C_n).$ By means of a coupling argument which
is also due to K.~Marton \cite{Mar96} (the so-called
\emph{Marton's coupling argument}), one can check that when $p=1,$
the product measure $\mu_1\otimes\cdots\otimes\mu_n$ satisfies
$T_1(C_1+\cdots+C_n),$ while when $p=2,$
$\mu_1\otimes\cdots\otimes\mu_n$ satisfies
$T_2(\max(C_1,\dots,C_n)).$ In particular, if $\mu$ satisfies
$T_1(C)$ then $\mu^{\otimes n}$ satisfies $T_1(nC).$ This
inequality deteriorates as $n$ grows. On the other hand, if $\mu$
satisfies $T_2(C)$ then $\mu^{\otimes n}$ also satisfies $T_2(C)$
and this still holds for the infinite product $\mu^{\otimes
\infty}.$
\\
By Jensen's inequality, we have $(\mathcal{T}_d)^2\leq
\mathcal{T}_{d^2}$ so that $T_2(C)$ implies $T_1(C).$ As the
standard normal law $\gamma$ satisfies $T_2(1),$ it is also shown
in \cite{Tal96a} that the standard normal law on $\mathbb{R}^n:$
$\gamma^{n},$ satisfies $T_2(1)$ and therefore $T_1(1)$ and the
concentration inequality
\begin{equation*}
\gamma^n(\{x; d(x,A)>r\})\leq
\exp\left[-\left(\frac{r}{\sqrt{2}}-\sqrt{\log 2}\right)^2\right]
\end{equation*}
for all measurable subset $A$ such that $\mu(A)\geq 1/2$ and all
$r\geq \sqrt{2\log 2}$ where $d$ is the Euclidean distance on
$\mathbb{R}^n.$ This concentration result holds for all $n$ and is
very close to the optimal concentration result obtained by means
of isoperimetric arguments (see M.~Ledoux's monograph
\cite{Led01}, Corollary 2.6) which is: $\gamma^n(\{x;
d(x,A)>r\})\leq e^{-r^2/2},$ for all $r\geq 0.$
\\
In view of (\ref{eq-10}) and of this optimal concentration
inequality, it now appears that with $\X= \mathbb{R}^n,$ $T_1(C)$
implies that $\mu$ concentrates at least as a normal law with
variance $C.$ One may say that $\mu$ performs a Gaussian
concentration when (\ref{eq-10}) holds for some $C.$

\par\medskip\noindent
 \textbf{Criteria for $T_1$.}\
It has recently been proved by H.~Djellout, A.~Guillin and L.~Wu
in  \cite{DGW03} that $\mu$ satisfies $T_1(C)$ for some $C$ if and
only if
\begin{equation}\label{eq-11}
    \IX e^{a_o d(x_o,x)^2}\,\mu(dx)<\infty
\end{equation}
 for some $a_o>0$ and some (and therefore all) $x_o$  in $\X.$ It
follows that \emph{(\ref{eq-11}) is a characterization of the
Gaussian concentration}. The proof of this result in \cite{DGW03}
relies on a dual characterization of $T_1$ which has been obtained
by S.~Bobkov and F.~G\"otze in \cite{BG99}. This characterization
is the following: $T_1(C)$ holds if and only if
\begin{equation}\label{eq-BG}
    \log \IX e^{s(\varphi-\langle\varphi,\mu\rangle)}\,d\mu
    \leq Cs^2/2,
\end{equation}
for all $s\geq 0$ and all bounded Lipschitz function $\varphi$
with $\|\varphi\|_\mathrm{Lip}\leq 1.$
\\
The criterion (\ref{eq-11}) has been recovered very recently by
F.~Bolley and C.~Villani in \cite{BV03} where the relation between
$C$ and $a_o$ is improved. This new proof relies on a
strengthening of CKP inequality  where weights are allowed in the
total variation norm. For a statement of this strengthened CKP
inequality, see Corollary \ref{res-32} below.

\subsection{Presentation of the results}

In this article, a larger class of transportation cost
inequalities is investigated. It appears that the transportation
cost inequalities $T_p$ defined by (\ref{eq-Tp}) enter the
following larger class of inequalities, which will also be called
transportation cost inequalities (TCIs):
\begin{equation}\label{eq-12}
\alpha(\mathcal{T}_c(\mu,\nu)) \leq H(\nu\mid\mu),\quad
\forall\nu\in\PX
\end{equation}
where $\alpha:[0,\infty)\rightarrow[0,\infty)$ is an
increasing\footnote{In the whole paper, by an increasing function
it is meant a nondecreasing function which may be constant on some
intervals.} function which vanishes at 0. The inequality
(\ref{eq-Tp}) corresponds $c=d^p$ with $\alpha(t)=t^{2/p}/(2C),$
$t\geq 0.$ Of course, one should rigorously restrict (\ref{eq-12})
to those $\nu\in\PX$ such that $\mathcal{T}_c(\mu,\nu)$ is
well-defined.

The aim of this paper is threefold.
\begin{itemize}
    \item[(i)]
    One proves TCIs by means of large deviation reasonings. The
    authors hope that this should provide a guideline for other functional inequalities.
    \item[(ii)]
    One obtains deviation results by means of TCIs.
    \item[(iii)]
    One extends already existing results, especially in the area
    of $T_1$-inequalities.
\end{itemize}

One says that we have a  $T_1$-inequality  if
\begin{equation*}
    \alpha(\mathcal{T}_d(\mu,\nu))\leq H(\nu\mid\mu),\ \forall \nu\in\mathcal{P}_d(\X).
   \tag{$T_1$}
\end{equation*}
where $d$ is a metric and $\mathcal{P}_d(\X)$ is the set of all
probability measures which integrate $d(x_o,x).$

\par\medskip\noindent
As regards item (i), it is no surprise that, because of the
relative entropy entering TCIs, Sanov theorem plays a crucial role
in our approach. Let
\[
L_n=\frac{1}{n}\sum_{i=1}^n\delta_{X_i}
\]
be the empirical
measure of an $n$-iid sample $(X_i)$ of the law $\mu\in\PX.$ Sanov
theorem states that the sequence $\{L_n\}_{n\geq 1}$ obeys the
large deviation principle with rate function $\nu\mapsto
H(\nu\mid\mu).$ The main idea is to control the deviations of the
nonnegative random variables $\mathcal{T}_c(\mu,L_n)$ as $n$ tends
to infinity. An easy heuristic description of this program is
displayed at Section \ref{sec-LD enter}. We obtain the

\begin{recipe}\label{recipe1}
Any increasing function $\alpha$ such that $\alpha(0)=0$ and
\begin{equation*}
        \logsup\P(\mathcal{T}_c(\mu,L_n)\geq t) \leq -\alpha(t)
\end{equation*}
for all $t\geq 0,$ satisfies the TCI (\ref{eq-12}).
\end{recipe}

Rigorously, one will have to require that $\alpha$ is a left
continuous function. This result will be proved at Theorem
\ref{res-abstract} and a weak version of it (with $\alpha$ convex)
is proved at Proposition \ref{res-51}.
\\
Not only TCIs can be derived with this recipe but also another
class of functional inequalities which we call Norm-Entropy
Inequalities (NEIs), see (\ref{eq-NEI}) for their definition. Let
us only emphasize in this introductory section that
$T_1$-inequalities are NEIs.

\par\medskip\noindent
As regards item (ii), \textbf{concentration inequalities} for
general measures and \textbf{deviation inequalities for empirical
processes} are derived by means of $T_1$-inequalities at Section
\ref{sec-applications}.

\par\medskip\noindent
As regards item (iii), the main technical (easy) result is Theorem
\ref{res-BG} which is an extension of Bobkov and G\"otze's
characterization of $T_1(C)$ stated at (\ref{eq-BG}). It gives a
\textbf{dual characterization} of all \emph{convex} TCIs: those
TCIs with $\alpha$ convex and increasing. Note that, up to the
knowledge of the authors, all known TCIs are convex. As a
consequence among others, one recovers the results of \cite{BV03}
about weighted CKP inequalities at Corollary \ref{res-32}.

\par\medskip\noindent
\textbf{Tensorization} of convex TCIs is also handled. The main
result on this topic is Theorem \ref{res-tensorization}. It states
that if $\alpha_1(\mathcal{T}_{c_1}(\mu_1,\nu_1))\leq
H(\nu_1\mid\mu_1)$ for all $\nu_1$ and
$\alpha_2(\mathcal{T}_{c_2}(\mu_2,\nu_2))\leq H(\nu_2\mid\mu_2)$
for all $\nu_2,$ then
$\alpha_1\square\alpha_2(\mathcal{T}_{c_1\oplus
c_2}(\mu_1\otimes\mu_2,\nu))\leq H(\nu\mid\mu_1\otimes\mu_2)$ for
all $\nu$ probability measure on the product space, where
$\alpha_1\square\alpha_2$ is the inf-convolution of $\alpha_1$ and
$\alpha_2.$

\par\medskip\noindent
\textbf{Integral criteria} are investigated in Section
\ref{sec-integral criteria}. It emerges from our analysis via
large deviations, that integral criteria only control the behavior
of $\alpha(t)$ in (\ref{eq-12}) for $t$ away from zero. As a
consequence, complete results are only derived for
$T_1$-inequalities. It is also proved that the function
$\alpha(t)$ of a $T_1$-inequality has a quadratic behavior for $t$
near zero. The integral criterion for $T_1$ is stated at Theorem
\ref{res-58}. It is the following:
\\
\textit{Let $d$ be a \lsc\ metric. Suppose that $a\geq 0$
satisfies $\IX e^{ad(x_o,x)}\,\mu(dx)\leq 2$ for some $x_o\in\X$
and that $\gamma$ is an increasing convex function which satisfies
$\gamma(0)=0$ and $\IX e^{\gamma(d(x_1,x))}\,\mu(dx)\leq B<\infty$
for some $x_1\in\X,$ then
\begin{equation*}
    \alpha(t)=\max\Big((\sqrt{at+1}-1)^2,2\gamma(t/2)-2\log
    B\Big),\ t\geq 0
\end{equation*}
satisfies} ($T_1$).
\\
Note that $(\sqrt{at+1}-1)^2=a^2 t^2/4+o_{t\rightarrow 0}(t^2)$ is
efficient for $t$ near zero, while $2\gamma(t/2)-2\log B$ is
efficient for $t$ away from zero.
\\
This theorem extends the integral criterion (\ref{eq-11}) of
\cite{DGW03} and \cite{BV03}.

\par\medskip\noindent
The last Section \ref{sec-abstract results} is devoted to abstract
results. In particular, the extended version Recipe \ref{recipe2}
of Recipe \ref{recipe1} is proved at Theorem \ref{res-abstract}.
The authors hope that the set of abstract results stated in this
section could be the starting point of the derivations of new
functional inequalities.

\tableofcontents

\section{Deriving $\mathcal{T}$-inequalities by means of large deviations. Heuristics}

The dual equality associated with the primal minimization problem
leading to $\mathcal{T}_c(\mu, \nu)$ is
\begin{equation}\label{eq-Kanto}
\mathcal{T}_c(\mu, \nu)=\sup_{(\psi,\varphi)\in\Phi_c}\left\{\IX
\psi\,d\mu+\IX \varphi\,d\nu\right\}
\end{equation}
where $\Phi_c$ is the set of all couples $(\psi,\varphi)$ of Borel
measurable bounded functions on $\X$ such that
$\psi(x)+\varphi(y)\leq c(x,y)$ for all $x,y\in\X.$ This result is
known as Kantorovich duality theorem and it holds true provided
that $c$ is \lsc. It still holds if $\Phi_c$ is replaced by
$C_b\cap \Phi_c$ which is the subset of all couples
$(\psi,\varphi)\in\Phi_c$ of continuous bounded functions. In the
special case where $c=d$ is a \lsc\ metric, the above dual
equality also holds with $\Phi_d$ the set of all couples
$(\psi,\varphi)$ of measurable (or continuous as well) bounded
functions such that $\psi=-\varphi$ and $\varphi$ is a
$d$-Lipschitz function with a Lipschitz constant less than 1. In
other words,
\begin{equation}\label{eq-KR}
\mathcal{T}_d(\mu, \nu)=\sup\left\{\IX \varphi\,d(\nu-\mu);
\varphi\in B(\X), \|\varphi\|_\mathrm{Lip}\leq
1\right\}:=\|\nu-\mu\|_{\mathrm{Lip}}^*
\end{equation}
where the space of all Borel measurable bounded functions on $\X$
is denoted $B(\X)$ and
$\|\varphi\|_{\mathrm{Lip}}=\sup_{x\not=y}\frac{|\varphi(x)-\varphi(y)|}{d(x,y)}$
is the usual Lipschitz seminorm. This result, known as
Kantorovich-Rubinstein's theorem, identifies the transportation
cost $\mathcal{T}_d(\mu,\nu)$ with the dual norm
$\|\nu-\mu\|_{\mathrm{Lip}}^*.$

\subsection{A larger class of transportation cost inequalities: $\mathcal{T}$-inequalities}
\label{sec-larger class of TCI}

After these considerations, it appears that the  transportation
cost inequality (\ref{eq-Tp}) enters the following larger class of
inequalities, which we call $\mathcal{T}$-inequalities:
\begin{equation}\label{eq-TIneq1}
\alpha(\mathcal{T}(\nu)) \leq H(\nu\mid\mu),\quad \forall\nu\in\N
\end{equation}
where $\alpha:[0,\infty)\rightarrow[0,\infty)$ is an increasing
function which vanishes at 0, $\N$ is a subset of $\PX$ and
$\mathcal{T}$ is defined by
\begin{equation}\label{def-T}
    \mathcal{T}(\nu)=\sup_{(\psi,\varphi)\in\Phi}\left\{\IX
\psi\,d\mu+\IX \varphi\,d\nu\right\}
\end{equation}
where $\Phi$ is a class of couples of functions $\pf$ with $\psi$
integrable with respect to $\mu$ and $\varphi$ integrable with
respect to $\nu.$ Note that (\ref{eq-TIneq1}) is a family of
inequalities where the value $+\infty$ is allowed with the
convention that
$\alpha(+\infty)=\lim_{t\rightarrow\infty}\alpha(t).$

 We are going to consider two cases which corresponds
to what will be called Transportation Cost Inequalities and
Norm-Entropy Inequalities.

 \textbf{Transportation Cost Inequalities.}\
We assume that  $c$ is a nonnegative \lsc\ cost function. The
space of all continous bounded functions on $\X$ is denoted
$C_b(\X).$ In the situation where $\Phi$ is equal to
\[
\Phi_c :=\{(\psi,\varphi)\in C_b(\X)\times C_b(\X);
\psi\oplus\varphi\leq c\}
\]
the family of inequalities (\ref{eq-TIneq1}) is called a
Transportation Cost Inequality (TCI). Indeed, the Kantorovich dual
equality (\ref{eq-Kanto}) states that
$$
\mathcal{T}(\nu)=\mathcal{T}_c(\mu,\nu)\in[0,\infty],
$$
for all $\nu\in \N\subset\PX.$ In this situation, inequality
(\ref{eq-TIneq1}) is
\begin{equation}\label{eq-TCI}
\alpha(\mathcal{T}_c(\mu,\nu)) \leq H(\nu\mid\mu),\quad
\forall\nu\in\N
\end{equation}

Suppose that there exists a nonnegative measurable function $\chi$
on $\X$ such that $c(x,y)\leq \chi(x)+\chi(y)$ for all $x,y\in\X$
and $\IX \chi\,d\mu<\infty.$ A natural set $\N$ is the set of all
probability measures $\nu$ such that $\IX\chi\,d\nu<\infty.$

\textbf{Norm-Entropy Inequalities.}\ Let $U$ be a set of
measurable functions on $\X$ such that $U=-U.$ Let us take
$\Phi=\Phi_U$ with
\[
\Phi_U :=\{(-\varphi,\varphi); \varphi\in U\}
\]
This gives
$$
\mathcal{T}(\nu)=\sup_{\varphi\in U}\IX
\varphi\,d(\nu-\mu):=\|\nu-\mu\|^*_U \in[0,\infty].
$$
In this case, inequality (\ref{eq-TIneq1}) is
\begin{equation}\label{eq-NEI}
\alpha(\|\nu-\mu\|^*_U) \leq H(\nu\mid\mu),\quad
\forall\nu\in\mathcal{P}_U
\end{equation}
where $\mathcal{P}_U$ is the set of all $\nu\in\PX$ such that $\IX
|\varphi|\,d\nu<\infty$ for all $\varphi\in U.$ The family of
inequalities (\ref{eq-NEI}) is called a Norm-Entropy Inequality
(NEI).

As a typical example, let $(F,\|\cdot\|)$ be a seminormed space of
measurable functions on $\X$ and $U:=\{\varphi\in F,
\|\varphi\|\leq 1\}$ its unit ball. Then, $\|\nu-\mu\|^*_U$ is the
dual norm of $\|\cdot\|.$

 \par\medskip
In the case where the cost function
of a TCI is a \lsc\ metric $d,$ the Kantorovich-Rubinstein theorem
(see (\ref{eq-KR})) states that
\[
\mathcal{T}_d(\mu,\nu)=\|\nu-\mu\|^*_{\mathrm{Lip}}
\]
for all $\mu, \nu\in\PX,$ where $\Phi_U$ is built with $F$ the
space all bounded $d$-Lipschitz functions on $\X$ endowed with the
seminorm $\|\cdot\|_{\mathrm{Lip}}.$ In this special important
case,  TCI and NEI match.

\subsection{Large deviations enter the game}\label{sec-LD enter}

At Sections \ref{sec-CvTI} and \ref{sec-abstract results},
$\mathcal{T}$-inequalities will be proved by means of a large
deviation approach. The integral functional $H(\cdot\mid\mu)$ will
be interpreted as the rate function of the large deviation
principle (LDP) of  the sequence
 of the empirical measures
 $$
 L_n=\frac{1}{n}\sum_{i=1}^n\delta_{X_i}
 $$
of an iid sample $(X_i)$ of the law $\mu$ ($\delta_x$ stands for
the Dirac measure at $x$). Indeed, by Sanov's theorem $\{L_n\}$
obeys the LDP in $\PX$ with the rate function
$$
I(\nu):=H(\nu\mid\mu),\ \nu\in \N.
$$

Roughly speaking, the sequence of random variables $\{L_n\}$ obeys
the LDP in $\N$ with the rate function $I$ if one has the
following collection of estimates
\begin{equation*}
   \P(L_n\in A) \asymp \exp[-n\inf_{\nu\in A} I(\nu)]
\end{equation*}
as $n$ tends to infinity, for any $A$  ``good" subset of $\N.$ Let
us introduce the nonnegative random variables
\begin{equation*}
    T_n=\mathcal{T}(L_n), \quad n\geq 1.
\end{equation*}
Suppose that $\mathcal{T}$ is regular enough for the sets
$A_t=\{\nu\in\N, \mathcal{T}(\nu)\geq t\},$ $t\geq 0,$ to be
``good" sets. This means that for all $t\geq 0,$
\[
\P(T_n\geq t)= \P(L_n\in A_t) \asymp \exp[-n i(t)]
\]
with $i(t)=\inf\{ I(\nu), \nu\in\N, \mathcal{T}(\nu)\geq
t\}\in[0,\infty].$ Suppose that $\alpha$ is a \emph{deviation
function} for the sequence $\{T_n\}$ in the sense that it is an
increasing nonnegative function on $[0,\infty)$ such that for all
$t\geq 0$
\begin{equation}\label{eq-deviation function}
    \logsup\P(T_n\geq t) \leq -\alpha(t).
\end{equation}

We obtain $\alpha(t)\leq i(t)$ for all $t$ and in particular with
$t=\mathcal{T}(\nu),$ we obtain for all $\nu\in\N,$
$\alpha(\mathcal{T}(\nu))\leq i(\mathcal{T}(\nu))\leq I(\nu).$
This is precisely the desired inequality (\ref{eq-TIneq1}).

The recipe is:
\begin{recipe}\label{recipe2}
Any deviation function $\alpha$ of $\{T_n\}$ satisfies the
$\mathcal{T}$-inequality (\ref{eq-TIneq1}).
\end{recipe}

Because of the sup entering the definition of $T_n=\sup_\Phi
(\langle \varphi, L_n\rangle+\langle\psi,\mu\rangle),$ one may
expect to get into troubles when trying to prove a full LDP for
$\{T_n\}.$ Fortunately, only the subclass of ``deviation sets"
$A_t=\{\nu\in\N, \mathcal{T}(\nu)\geq t\},$ $t\geq 0,$ will be
really useful.

This line of reasoning will be put on a solid ground at Theorem
\ref{res-BG}, Proposition \ref{res-51} and Theorem
\ref{res-abstract}.

\subsection{An example: CKP inequality}

As a simple illustration, we propose to prove CKP inequality by
searching a deviation function $\alpha$ in the sense of
(\ref{eq-deviation function}). This is not intended to be the
shortest proof, but only an illustration of the proposed method.
Recall that CKP inequality is
\begin{equation}\label{eq-Pinsker}
    \frac{1}{2}\|\nu-\mu\|_{\mathrm{TV}}^2\leq H(\nu\mid\mu),\
    \forall \nu\in\PX
\end{equation}
where $\|\xi\|_{\mathrm{TV}}$ is the total variation of the signed
bounded measure $\xi.$ As
\begin{equation*}
  \|\xi\|_{\mathrm{TV}}
  = \sup\left\{\IX \varphi\,d\xi,; \varphi \textrm{ measurable such that }
  \|\varphi\|:=\sup_{x\in\X}|\varphi(x)|\leq
   1\right\},
\end{equation*}
(\ref{eq-Pinsker}) is the NEI with $F=B(\X)$ the space of bounded
measurable functions furnished with the uniform norm
$\|\varphi\|:=\sup_{x\in\X}|\varphi(x)|,$ $\N=\PX$ and
$\alpha(t)=t^2/2.$

Consider an iid sample $(X_i)$ of the law $\mu$ and its associated
sequence of empirical measures
$L_n=\frac{1}{n}\sum_{i=1}^n\delta_{X_i}.$ For all $n$ and all
 $\varphi\in U=\{\varphi \in B(\X);\|\varphi\|\leq 1\},$ define
 the random variable
 \[
T^\varphi_n=\langle \varphi,
L_n-\mu\rangle=\frac{1}{n}\sum_{i=1}^n Y_i^\varphi
 \]
where $Y_i^\varphi=\varphi(X_i)-\E \varphi(X_i).$  Cram\'er's
theorem states that $\{T^\varphi_n\}$ obeys the LDP in $\R$ with
 rate function $\Lambda^*_\varphi:$  the convex conjugate of the
log-Laplace transform $\Lambda_\varphi(s)=\log\E e^{sY^\varphi},$
$s\in\R.$ Recall that the convex conjugate of $f$ is defined by
$f^*(t)=\sup_{s\in\R}\{st-f(s)\}\in (-\infty,\infty],$ $t\in\R.$

Sanov's theorem holds in $\PX$ with the weak topology
$\sigma(\PX,B(\X)).$ As, $\nu\in\PX\mapsto
\langle\varphi,\nu-\mu\rangle$ is $\sigma(\PX,B(\X))$-continuous
for all $\varphi\in B(\X),$ one can apply the contraction
principle. It gives us for all $t$
\begin{equation*}
    \Lambda_\varphi^*(t)=\inf\{H(\nu\mid\mu); \nu\in\PX:
    \langle\varphi,\nu-\mu\rangle=t\},
\end{equation*}
which in turn implies that for all $\varphi\in B(\X),$
\begin{equation*}
    \Lambda_\varphi^*(\langle\varphi,\nu-\mu\rangle)\leq H(\nu\mid\mu), \ \forall
    \nu\in\PX.
\end{equation*}

As $Y^\varphi$ takes its values in $[\E Y^\varphi-1,\E
Y^\varphi+1],$ by Hoeffding's inequality we have
\begin{equation}\label{eq-Hoeffding}
    \Lambda_\varphi(s)\leq s^2/2
\end{equation}
for all real $s.$ It follows that $\Lambda^*_\varphi(t)\geq
\sup_{s\in\R}\{st-s^2/2\}=t^2/2$ for all real $t.$ Hence, we have
proved that for all $\varphi\in U,$
\begin{equation*}
 \alpha(\langle\varphi,\nu-\mu\rangle)\leq H(\nu\mid\mu), \ \forall
    \nu\in\PX
\end{equation*}
with $\alpha(t)=t^2/2.$ It follows that $\alpha(\sup_{\varphi\in
U}\langle\varphi,\nu-\mu\rangle)\leq H(\nu\mid\mu)$ for all
$\nu\in\PX,$ which is CKP inequality (\ref{eq-Pinsker}).

\textbf{Some comments.}\ In this proof, something interesting
occured. Let us denote $T_n:=\sup_\varphi T_n^\varphi,$
$\beta(t)=-\logsup\P(T_n\geq t)$ and
$J_\varphi(t)=-\logsup\P(T_n^\varphi\geq t)$ the deviation
functions of $T_n$ and $T_n^\varphi.$ As $T_n\geq T_n^\varphi$ for
all $\varphi,$ we have $\beta\leq \inf_\varphi J_\varphi.$ This
means that a priori $\inf_\varphi J_\varphi$ could be too
 large to be the $\alpha$ of the NEI.
On the other hand, by (\ref{eq-Hoeffding}):
$\sup_\varphi\Lambda_\varphi(s)\leq \Lambda(s):=s^2/2$ for all
$s\geq 0,$ so that $t^2/2=\Lambda^*(t)\leq \inf_\varphi
J_\varphi(t).$
\\
Nevertheless, we have shown that $\Lambda^*$ is a convenient
function $\alpha$ for our NEI.
\\
It will shown in a more general setting, at Theorem
\ref{res-reassuring}, that the convex \lsc\ envelope of
$\inf_\varphi J_\varphi$ is the best \emph{increasing convex}
function $\alpha$ for this NEI.

\section{Convex $\mathcal{T}$-inequalities. A dual characterization}\label{sec-CvTI}

In the rest of the paper (except Section \ref{sec-abstract
results}) our attention is restricted to those
$\mathcal{T}$-inequalities (\ref{eq-TIneq1}) where the function
$\alpha$ is increasing and convex. In this case, (\ref{eq-TIneq1})
is said to be a convex $\mathcal{T}$-inequality.

\subsection{Sanov's theorem}\label{sec-Sanov}

This theorem will be central for the proof of the main result of
this section which is stated at Theorem \ref{res-BG}.
\\
Let the probability measure $\mu$ on $\X$ be given. We consider a
sequence of independent $\X$-valued random variables $(X_i)_{i\geq
1}$ identically distributed with law $\mu.$ For any $n$ the
empirical measure of this sample is
$$
L_n=\frac{1}{n}\sum_{i=1}^n\delta_{X_i}\in\PX.
$$

We introduce the function space
\begin{equation}\label{eq-34}
    \mathcal{F}_{\mathrm{exp}}(\mu)=\left\{\varphi:\X\rightarrow\R;
\varphi \textrm{ measurable, } \IX \exp(a |\varphi|)\,d\mu<\infty
 \textrm{ for all } a>0\right\}
\end{equation}

of all the functions which admit exponential moments of all orders
with respect to the measure $\mu.$ We denote
\[
\N_{\mathrm{exp}}(\mu)=\left\{\nu\in\PX; \IX
|\varphi|\,d\nu<\infty \textrm{ for all }
\varphi\in\mathcal{F}_{\mathrm{exp}}(\mu)\right\}
\]
the set of all probability measures which integrate every function
of $\mathcal{F}_{\mathrm{exp}}(\mu).$

The set $\PX$ is furnished with the cylinder $\sigma$-field
generated by the functions $\nu\mapsto \langle
\varphi,\nu\rangle,$ $\varphi\in\mathcal{F}_{\mathrm{exp}}(\mu).$

\begin{thm}[A version of Sanov's theorem]\label{res-Sanov}
The effective domain of $H(\cdot\mid\mu)$ is included in
$\N_{\mathrm{exp}}(\mu)$ and the sequence $\{L_n\}$ obeys the
large deviation principle with rate function $H(\cdot\mid\mu)$ in
$\N_{\mathrm{exp}}(\mu)$ equipped with the weak topology
$\sigma(\N_{\mathrm{exp}}(\mu),\mathcal{F}_{\mathrm{exp}}(\mu)).$

This means that for all measurable subset $A$ of
$\N_{\mathrm{exp}}(\mu),$  we have
\begin{eqnarray*}
  \loginf \P(L_n\in A) &\geq& -\inf_{\nu\in \inter A}H(\nu\mid\mu) \textrm{\quad and} \\
  \logsup \P(L_n\in A) &\leq& -\inf_{\nu\in \cl A}H(\nu\mid\mu)
\end{eqnarray*}
where $\inter A$ and $\cl A$ are the interior and closure of $A.$
\end{thm}

\proof The proof is a variation of the classical proof of Sanov's
theorem based on projective limits of LD systems (see \cite{DZ},
Thm 6.2.10). For two distinct detailed proofs of the present
theorem, see (\cite{EiS}, Theorem 1.7) or (\cite{LeoN02},
Corollary 3.3).
\endproof

\subsection{The class of functions $\mathcal{C}$}

The functions $\alpha$ to be considered are assumed to be convex.
Since $\alpha $ is also left continuous and increasing, we
consider the following class of functions.

\begin{definition}[of $\mathcal{C}$]
The class $\mathcal{C}$ consists of all the functions $\alpha$ on
$[0,\infty)$ which are \emph{convex} increasing, left continuous
with $\alpha(0)=0.$
\end{definition}

For any $\alpha$ belonging to the class $\mathcal{C},$ denoting
$t_*=\sup\{t\geq 0; \alpha(t)<\infty\},$ $\alpha$ is continuous on
$[0,t_*)$ and $\lim_{t\uparrow t_*}\alpha(t)=\alpha(t_*).$

The convex conjugate of a function $\alpha\in\mathcal{C}$ is
replaced by the monotone conjugate $\alpha^{\mc}$ defined by
\[
\alpha^{\mc}(s)=\sup_{t\geq 0}\{st-\alpha(t)\}, s\geq 0
\]
where the supremum in taken on $t\geq 0$ instead of $t\in\R.$ In
fact, if $\alpha$ is extended by
$\widetilde{\alpha}(t)=\left\{\begin{array}{ll}
  \alpha(t) & \textrm{if } t\geq 0 \\
  0 & \textrm{if } t\leq 0 \\
\end{array}\right.$ then the usual convex conjugate of
$\widetilde{\alpha}$ is
$\widetilde{\alpha}^*(s)=\left\{\begin{array}{ll}
  \alpha^{\mc}(s) & \textrm{if } s\geq 0 \\
  +\infty & \textrm{if } s< 0 \\
\end{array}\right..$
As $\widetilde{\alpha}$ is convex and \lsc, we have
$\widetilde{\alpha}^{**}=\widetilde{\alpha}.$ From this, it is not
hard to deduce the following result.

\begin{prop}\label{res-class C}
For any function $\alpha$ on $[0,\infty),$ we have
\begin{itemize}
    \item[(a)]
    $\alpha\in\mathcal{C}\Leftrightarrow\alpha^{\mc}\in\mathcal{C}$
    \item[(b)]
    $\alpha\in\mathcal{C}\Rightarrow\alpha^{\mc\mc}=\alpha.$
\end{itemize}
\end{prop}

\subsection{A convex criterion}

Theorem \ref{res-BG} below is a criterion for a convex
$\mathcal{T}$-inequality to hold. It extends two well-known results
of S.~Bobkov and F.~G\"otze (\cite{BG99}, Theorem 1.3 and statement (1.7)).

Let $\mathcal{F}$ be a vector space of measurable functions
$\varphi$ on $\X$ such that
\begin{equation}\label{hyp-F}
    \IX e^{\varphi}\,d\mu<\infty,\quad  \forall
    \varphi\in \mathcal{F}.
\end{equation}

Let $\NF$ be the set of all probability measures which integrate
$\mathcal{F}:$
\begin{equation*}
    \NF=\left\{\nu\in\PX; \IX |\varphi|\,d\nu<\infty,\
    \forall\varphi\in\mathcal{F}\right\}.
\end{equation*}

Clearly, if the class $\Phi$ entering the definition of
$\mathcal{T}(\nu)$ satisfies
\begin{equation}\label{hyp-Phi}
   (0,0)\in \Phi\subset \mathcal{F}\times\mathcal{F},
\end{equation}
the function $\mathcal{T}$ is a well defined $[0,\infty]$-valued
function on $\NF.$

Let $\Lambda_\phi(s)$ be the log-Laplace transform of
$\varphi(X)+\E\psi(X)$ where $X$ admits $\mu$ as its law. We have
for all real $s,$
\[
\Lambda_\phi(s)=\log\IX
\exp[s(\varphi(x)+\langle\psi,\mu\rangle)]\,\mu(dx)
\]

\begin{thm}\label{res-BG}
We assume (\ref{hyp-F}) and (\ref{hyp-Phi}). Let us consider the
following statements where $\alpha$ is any function  in
$\mathcal{C}:$
\begin{enumerate}
    \item[(a)]
    $\alpha(\mathcal{T}(\nu))\leq H(\nu\mid\mu),$\quad
    $\forall\nu\in\NF.$
    \item[(b)]
        $\Lambda_\phi(s)\leq
    \alpha^{\mc}(s),$\quad
    $\forall s\geq 0,$ $\forall\phi\in\Phi.$
    \item[(c)]
    $\alpha(t)\leq \Lambda_\phi^*(t), \forall t\geq 0, \forall
    \phi\in\Phi.$
    \item[(d)]
    $\logsup\P(\langle\varphi,L_n \rangle+\langle\psi,\mu\rangle \geq  t)
    \leq -\alpha(t),$\quad
    $\forall t\geq 0,$ $\forall(\psi,\varphi)\in\Phi.$
    \item[(e)]
    $\forall n\geq 1,$\
    $\frac{1}{n}\log\P(\langle\varphi,L_n \rangle+\langle\psi,\mu\rangle \geq  t)
    \leq -\alpha(t),$\quad
    $\forall t\geq 0,$ $\forall(\psi,\varphi)\in\Phi.$ \\
\end{enumerate}
Then, we have
$
(a)\Leftrightarrow (b)\Leftrightarrow (c)
$
and $(e)\Rightarrow (d)\Rightarrow (a).$
\\
If it is assumed in addition that for all $\pf\in\Phi,$
\begin{equation}\label{eq-37}
    \IX (\varphi(x)+\psi(x))\,\mu(dx)\leq 0
\end{equation}

then, we have
$
(a)\Leftrightarrow (b)\Leftrightarrow (c)\Leftrightarrow
(d)\Leftrightarrow (e).
$
\end{thm}

The most useful statement of this theorem is the criterion
$(b)\Rightarrow (a).$
\\
Clearly, the requirement (\ref{eq-37}) holds for all NEIs. It also
holds for TCIs under the assumption  that $c$ satisfies
\begin{equation}\label{eq-21}
    c(x,x)=0,\ \forall x\in\X.
\end{equation}
When working with TCIs, this will be assumed in the sequel.

\proof

Possibly considering the vector space $\mathcal{F}'$ spanned by
$\mathcal{F}\cup C_b(\X)$ instead of $\mathcal{F},$ one can assume
that $\mathcal{F}$ separates $\NF.$ Indeed, the assumptions
(\ref{hyp-F}) and (\ref{hyp-Phi}) still hold with $\mathcal{F}'$
instead of $\mathcal{F}$ and we clearly have
$\mathcal{P}_{\mathcal{F}'}=\NF.$ Hence, we assume without loss of
generality that $\mathcal{F}$ separates $\NF.$ As a consequence,
the weak topology $\sigma(\NF,\mathcal{F})$ is Hausdorff: this is
necessary to derive LDPs away from compactness troubles.

Note that  the assumption (\ref{hyp-F}) is equivalent to
$\mathcal{F}\subset\mathcal{F}_{\mathrm{exp}}(\mu).$ It follows
that under this assumption, Sanov's Theorem \ref{res-Sanov}
implies that $\{L_n\}$ obeys the LDP in $\NF$ equipped with
$\sigma(\NF,\mathcal{F})$ with $H(\cdot\mid\mu)$ as its rate
function.

Consider, for any $(\psi,\varphi):=\phi\in\Phi$ and $n\geq 1,$
\begin{equation}\label{eq-33}
    T_n^{\phi}=\langle\varphi,L_n\rangle+\langle\psi,\mu\rangle
    =\frac{1}{n}\sum_{i=1}^n (\varphi(X_i)+\E\psi(X_i))
\end{equation}
so that $T_n:=\mathcal{T}(L_n)=\sup_{\phi\in\Phi}T_n^{\phi}.$
Cram\'er's theorem states that $\{T_n^\phi\}$ obeys the LDP in
$\R$ with
\[
\Lambda_\phi^*(t)=\sup_{s\in\R}\{st-\Lambda_\phi(s)\},\ t\in\R
\]
as its rate function. In particular, for all real $t$
\begin{eqnarray}\label{eq-32}
\nonumber   -\inf_{u> t}\Lambda_\phi^*(u)&\leq&\loginf\P(T_n^\phi> t) \\
 &\leq& \logsup\P(T_n^\phi\geq t)\leq -\inf_{u\geq t}\Lambda_\phi^*(u)
\end{eqnarray}

Because of assumption (\ref{hyp-Phi}), the mapping
$f_\phi:\nu\in\NF\mapsto\langle\varphi,\nu\rangle+\langle\psi,\mu\rangle\in\R$
is continuous for every $\pf\in\Phi.$ As $T_n^\phi=f_\phi(L_n),$
one can apply the contraction principle which gives us for all
real $t$
\begin{equation}\label{eq-31}
    \Lambda_\phi^*(t)=\inf\{H(\nu\mid\mu); \nu\in\NF:
\langle\varphi,\nu\rangle+\langle\psi,\mu\rangle=t\}.
\end{equation}

 \par\medskip\noindent $[(a)\Leftrightarrow (c)]:$\quad
\begin{eqnarray*}
 (a) &\stackrel{(i)}{\Leftrightarrow}& \alpha\left(\sup_\phi (\langle\varphi,\nu\rangle+\langle\psi,\mu\rangle)\right)\leq H(\nu\mid\mu), \forall \nu\in\NF \\
   &\stackrel{(ii)}{\Leftrightarrow}&  \alpha(\langle\varphi,\nu\rangle+\langle\psi,\mu\rangle)\leq H(\nu\mid\mu), \forall \nu\in\NF ,\forall \phi\in\Phi\\
   &\Leftrightarrow& \alpha(t)\leq H(\nu\mid\mu), \forall t\in\R,\forall \phi\in\Phi,  \forall \nu\in\NF: \langle\varphi,\nu\rangle+\langle\psi,\mu\rangle=t\\
   &\Leftrightarrow& \alpha(t)\leq \inf\{H(\nu\mid\mu); \nu\in\NF: \langle\varphi,\nu\rangle+\langle\psi,\mu\rangle=t\}, \forall t\in\R,\forall \phi\in\Phi  \\
   &\stackrel{(iii)}{\Leftrightarrow}& \alpha\leq \Lambda_\phi^*\\
   &\Leftrightarrow& (c)
\end{eqnarray*}
The equivalence (i) follows from the definition (\ref{def-T}) of
$\mathcal{T},$ (ii) holds true because $\alpha$ is increasing and
left continuous while (iii) follows from (\ref{eq-31}).

\par\medskip\noindent $[(b)\Leftrightarrow (c)].$\quad
In order to work with usual convex conjugates instead of monotone
conjugates, let us take $\alpha^*(s)=+\infty$ for all $s<0.$ It
follows that $\alpha$ is extended by $\alpha(t)=0,$ for all $t\leq
0$ and $\alpha^\mc(s)=\alpha^*(s)$ for all $s\geq 0.$

Let us prove $(c)\Rightarrow (b).$ With the above convention,
statement (c) is equivalent to
\begin{equation}\label{eq-38}
\alpha(t)\leq \Lambda_\phi^*(t), \forall t\in\R, \forall
    \phi\in\Phi.
\end{equation}

As, $\Lambda_\phi$ is  convex and \lsc, we have:
$\Lambda_\phi^{**}=\Lambda_\phi.$ Hence, taking the convex
conjugates on both sides of (\ref{eq-38}) one obtains that
$\Lambda_\phi\leq \alpha^*$ which entails (b).

Let us prove $(b)\Rightarrow (c).$
 As $\alpha$ is in $\mathcal{C},$ its extension (still denoted
by $\alpha$) is convex and \lsc, so that $\alpha^{**}=\alpha.$
Therefore, taking the conjugate of (b) leads to
$\alpha\leq\Lambda_\phi^*$ which is (c).

The convexity of $\alpha$ has been used to obtain $ (b)\Rightarrow
(c) $ and it won't be used anywhere else.

\par\medskip\noindent $[(e)\Rightarrow (d)\Rightarrow (a)].$\quad
 As
$(e)\Rightarrow (d)$ is obvious and $(a)\Leftrightarrow (c),$ all
we have to show is $(d)\Rightarrow(c).$
\\
 Let $m=\E Y=\langle
\varphi+\psi,\mu\rangle.$ For all $t\leq m,$ we have
 $\inf_{u>t}\Lambda_\phi^*(u)=\inf_{u\geq t}\Lambda_\phi^*(u)=0.$
As $\Lambda_\phi^*$ is convex, it is continuous on $(t_-,t_+)$ the
interior of its effective domain. Therefore, we have for all
$t\not =t_+,$  $\inf_{u>t}\Lambda_\phi^*(u)=\inf_{u\geq
t}\Lambda_\phi^*(u).$ Together with (\ref{eq-32}), this gives for
all $t\not =t_+,$
\begin{equation*}
    -\lim_{n\rightarrow\infty}\frac{1}{n}\log\P(T_n^\phi\geq t)
  = \inf_{u>t}\Lambda_\phi^*(u)=\inf_{u\geq t}\Lambda_\phi^*(u)\\
  = \left\{%
\begin{array}{ll}
    0, & \hbox{if } t\leq m \\
    \Lambda_\phi^*(t), & \hbox{if } t\geq m \\
\end{array}%
\right. =\Lambda_\phi^{\mc}(t).
\end{equation*}
Consequently, considering $\Gamma(t)=\Lambda_\phi^{\mc}(t)$ if
$t\not =t_+$ and $\Gamma(t_+)=+\infty$ (if $t_+<\infty$), we have
\begin{eqnarray*}
  (d) &\Rightarrow& \alpha(t)\leq \Lambda_\phi^{\mc}(t),\ \forall t\not =t_+ \\
   &\Rightarrow& \alpha\leq \Gamma \\
    &\Rightarrow& \mathrm{ls\ }\alpha\leq \mathrm{ls\ }\Gamma \\
   &\Rightarrow&  \alpha\leq \Lambda_\phi^{\mc}
\end{eqnarray*}
where $\mathrm{ls\ }\alpha$ and $\mathrm{ls\ }\Gamma$ are the
\lsc\ envelopes of  $\alpha$ and $\Gamma,$ and the last
implication holds since $\alpha$ is \lsc\ and $\mathrm{ls\
}\Gamma=\Lambda_\phi^{\mc}.$ As
$\Lambda_\phi^{\mc}\leq\Lambda_\phi^*,$ we have the desired
result.

\par\medskip\noindent $[(a)\Leftrightarrow (b)\Leftrightarrow (c)\Leftrightarrow (d)\Leftrightarrow (e)].$\quad
Let us assume (\ref{eq-37}). To obtain the stated series of
equivalences, it remains to prove $(c)\Rightarrow (e).$
\\
By (\ref{eq-33}), $T_n^{\phi} =\frac{1}{n}\sum_{i=1}^n Y_i$ with
$Y_i=\varphi(X_i)+\E\psi(X_i).$ The standard proof of the upper
bound of Cram\'er's theorem is based on an optimization of a
collection of exponential Markov inequalities, as follows. For all
real $t,$ all $n$ and all $s\geq 0,$
\begin{eqnarray*}
  \P\left(\frac{1}{n}\sum_{i=1}^n Y_i\geq t\right) &\leq& \P\left(\exp[s\sum_{i=1}^nY_i]\geq e^{nst}\right) \\
   &\leq& e^{-nst}\E \exp [s\sum_{i=1}^n Y_i]\\
   &=& \exp[n(\Lambda_\phi(s)-st)]
\end{eqnarray*}
Optimizing on $s\geq 0,$ one obtains that
\[
\frac{1}{n}\log\P(T_n^\phi\geq t)\leq -\Lambda_\phi^{\mc}(t),\
\forall t\in\R, \forall \phi\in\Phi, \forall n\geq 1.
\]
But, assumption (\ref{eq-37}) implies that $m\leq 0$ so that
$\Lambda_\phi^\mc(t)=\Lambda_\phi^*(t)$ for all $t\geq 0.$ It
follows immediately that $(c)\Rightarrow (e).$ This completes the
proof of the theorem.
\endproof

\subsection{Convex Transportation Cost Inequalities}

In the special case of TCIs, we have $\Phi=\Phi_c=\{\pf;\psi,
\varphi\in C_b(\X): \psi\oplus\varphi\leq c\}.$ Optimal
transportation theory (see \cite{Vill}) indicates that $\Phi_c$
may be replaced with the smaller sets
$\{(-\varphi,Q^c\varphi);\varphi\in C_b(\X)\}$ or
 $\{(-\varphi,Q^c\varphi);\varphi \textrm{ \lsc\ and bounded on }\X\}$
where
\[
Q^c\varphi(y)=\inf_{x\in\X}\{\varphi(x)+c(x,y)\},\ y\in\X
\]
without any change in the value of $\mathcal{T}_c.$ One easily
proves that if (\ref{eq-21}) is satisfied:  $c(x,x)=0$ for all
$x\in\X,$ then $\sup|Q^c\varphi|\leq \sup|\varphi|.$ If $c$ is
continuous, then $Q^c\varphi$ is measurable as an upper
semicontinuous function. If $c$ is only assumed to be \lsc,
$Q^c\varphi$ is still measurable if $\varphi$ is \lsc\ and bounded
(but the proof of this result is technical). Anyway,
$Q^c\varphi\in B(\X)$ (is a bounded measurable function) as soon
as $\varphi$ is \lsc\ and bounded. In particular, assumptions
(\ref{hyp-F}) and (\ref{hyp-Phi}) hold with $\mathcal{F}=B(\X).$

Now, as a corollary of Theorem \ref{res-BG}, we have the following
result.

\begin{cor}\label{res-34}
 Whenever $\alpha\in \mathcal{C},$ the transportation cost inequality (\ref{eq-TCI})
holds in $\N=\PX$ if and only if
\begin{equation*}
  \log\IX
e^{s[Q^c\varphi(y)-\langle\varphi,\mu\rangle]}\,\mu(dy) \leq
\alpha^{\mc}(s)
\end{equation*}
for all $s\geq 0$ and all $\varphi\in C_b(\X).$

If in addition $c$ is continuous, the same result holds when
$\varphi\in C_b(\X)$ is replaced with $\varphi\in B(\X):$ the set
of all measurable bounded functions on $\X.$
\end{cor}

\subsection{Convex Norm-Entropy inequalities}

In the special case of NEIs, we have $\Phi=\{(-\varphi,\varphi);
\varphi\in U\}$ and  Theorem \ref{res-BG} specializes as
follows.

\begin{thm}\label{res-BG-NEI} Suppose that $U$ satisfies
\[
\IX e^{a|\varphi|}\,d\mu<\infty, \forall \varphi\in U, \forall
a>0.
\]
Let $\alpha$ be in $\mathcal{C}.$  Then, the norm-entropy
inequality (\ref{eq-NEI})
\begin{equation*}
\alpha(\|\nu-\mu\|^*_U) \leq H(\nu\mid\mu),\quad \forall\nu\in
\mathcal{P}_U
\end{equation*}
 holds if and only if
\begin{equation}\label{eq-35}
    \Lambda_\varphi(s):=\log\IX
e^{s[\varphi(x)-\langle\varphi,\mu\rangle]}\,\mu(dx) \leq
\alpha^{\mc}(s)
\end{equation}
for all $s\geq 0$ and all $\varphi\in U.$
\end{thm}

\par\medskip
Specializing  Theorem \ref{res-BG-NEI} by taking $U$ to be the set
of all 1-Lipschitz measurable bounded functions with respect some
measurable metric $d,$ one obtains the following characterization
of convex $T_1$-inequalities.

\begin{thm}[$T_1$-inequality]\label{res-BG-T1} Let $d$ be a \lsc\ metric on $\X$
such that
\[
\IX e^{a_o d(x_o,x)}\,\mu(dx)<\infty,
\]
for some $a_o>0$  and some (and therefore all)  $x_o\in\X.$ Let
$\alpha$ be in $\mathcal{C}.$  Then,
\begin{equation*}
\alpha(\mathcal{T}_d(\mu,\nu)) \leq H(\nu\mid\mu),
\end{equation*}

for all $\nu\in \PX $ such that $\IX d(x_o,x)\,\nu(dx)<\infty$ if
and only if
\begin{equation}\label{eq-35bis}
    \Lambda_\varphi(s):=\log\IX
e^{s[\varphi(x)-\langle\varphi,\mu\rangle]}\,\mu(dx) \leq
\alpha^{\mc}(s)
\end{equation}
for all $s\geq 0$ and all measurable bounded Lipschitz function
$\varphi$ such that $\|\varphi\|_\mathrm{Lip}\leq 1.$
\end{thm}

The following simple result asserts that the functions $\alpha$ of
NEIs cannot grow faster than $at^2$ for $t$ near zero.

\begin{prop}\label{res-35}
Assuming that $F$ contains functions which are not $\mu$-a.e.
constant, the function $\alpha$ of a convex norm-entropy
inequality (\ref{eq-NEI}) satisfies
\begin{equation}\label{eq-39}
0\leq\alpha(t)\leq at^2, \forall 0\leq t\leq t_1
\end{equation}
for some $a>0$ and $t_1>0.$
\end{prop}

\proof Let $\varphi_o$ be a non constant function in $U.$ Then,
$\sigma_o^2:=\IX (\varphi(x)-\langle\varphi,\mu\rangle)^2\,d\mu>0$
and for any $0<\sigma_1^2<\sigma_o^2$ there exists $s_1>0$ such
that $\Lambda_{\varphi_o}(s)=\sigma_o^2s^2/2+o(s^2)\geq
\sigma_1^2s^2/2,$ for all $0\leq s\leq s_1.$ Let $\theta_1(s)$
match with $\sigma_1^2s^2/2$ on $[0,s_1]$ and be extended on
$[s_1,\infty)$ by the tangent affine function of
$s\mapsto\sigma_1^2s^2/2$ at $s=s_1.$ As $\Lambda_{\varphi_o}$ is
convex, we have $\theta_1(s)\leq \Lambda_{\varphi_o}(s)$ for all
$s\geq 0.$

Together with (\ref{eq-35}), we obtain $\theta_1\leq
\alpha^{\mc}.$ Taking the monotone conjugates on both sides of
this inequality provides us with
\begin{equation*}
    \alpha(t)\leq \theta_1^{\mc}(t)=\left\{%
\begin{array}{ll}
   t^2/(2 \sigma_1^2), & \hbox{if }0\leq t \leq s_1\sigma_1^2 \\
    +\infty, & \hbox{if } t > s_1\sigma_1^2\\
\end{array}%
\right.
\end{equation*}
from which the desired result follows.
\endproof

To explore some  consequences of Theorem \ref{res-BG-NEI} (see
Corollaries \ref{res-31} and \ref{res-32} below) one needs the
notion of Orlicz space associated with the exponential function.
It appears that the space $\mathcal{F}_{\mathrm{exp}}(\mu)$
introduced at (\ref{eq-34}) is the Orlicz space
$$
\left\{\varphi:\X\rightarrow\R; \textrm{ measurable,
}\IX\rho(a\varphi)\,d\mu<\infty \textrm{ for all } a>0\right\}
$$
where $\mu$-almost equal functions are not identified and $\rho$
is the Young function
\[
\rho(s)=e^{|s|}-1,\quad s\in\R.
\]
Its Orlicz norm is defined by
\begin{eqnarray}\label{eq-Orlicz norm}
   \|\varphi\|_\rho &:=& \inf\left\{b>0;
\IX\rho\left(\frac{\varphi}{b}\right)\,d\mu\leq 1\right\} \\
  \nonumber  &=& \inf\left\{b>0; \IX e^{|\varphi|/b}\,d\mu\leq 2\right\}
\end{eqnarray}

and considering the usual dual bracket
$\langle\eta,\varphi\rangle=\IX \eta\varphi\,d\mu,$ its
topological dual space is isomorphic to
\begin{eqnarray*}
  L_{\rho^*}(\mu) &=& \left\{\eta:\X\rightarrow\R; \textrm{ measurable, }\IX\rho^*(a\eta)\,d\mu<\infty \textrm{ for some
}a>0\right\} \\
   &=& \left\{\eta:\X\rightarrow\R; \textrm{ measurable, }\IX|\eta|\log|\eta|\,d\mu<\infty\right\}
\end{eqnarray*}

where $\rho^*$ is the convex conjugate of $\rho:$
\[
\rho^*(t)=\left\{%
\begin{array}{ll}
    |t|\log |t| -|t|+1, & \hbox{if }|t|\geq 1 \\
    0, & \hbox{if } |t|\leq 1\\
\end{array}%
\right.
\]
and $\mu$-almost equal functions are identified. Note that the
effective domain of $H(\cdot\mid\mu)$ is included in the set of
all probability measures $\nu$ which are absolutely continuous
with respect to $\mu$ and such that $\frac{d\nu}{d\mu}\in
L_{\rho^*}(\mu) .$

Let us state a  useful technical lemma, which will play a role
that is similar to the role that Hoeffding's inequality
(\ref{eq-Hoeffding}) played during the proof of CKP inequality.

\begin{lemma}[A Bernstein type inequality]\label{res-Bernstein-Orlicz}
For any measurable function $\varphi$ such that $\IX
e^{a_o|\varphi|}\,d\mu<\infty$ for some  $a_o>0,$ we have
$\|\varphi\|_\rho<\infty$ and
\[
\Lambda_\varphi(s)\leq \frac{\|\varphi\|_\rho^2\,
s^2}{1-\|\varphi\|_\rho\, s},\quad\forall\ 0\leq s<
1/\|\varphi\|_\rho.
\]
It follows that, if $U$ is a uniformfy $\|\cdot\|_\rho$-bounded
set of functions: $\sup_{\varphi\in U}\|\varphi\|_\rho\leq
M<\infty,$ then
\[
\Lambda_\varphi(s)\leq \frac{M^2 s^2}{1-M s},\quad\forall\ 0\leq
s< 1/M, \forall \varphi\in U.
\]
\end{lemma}

\proof By the definition of $\beta:=\|\varphi\|_\rho,$ we have
$1\geq\IX \rho(\varphi/\beta)\,d\mu=\sum_{k\geq 1} \langle
|\varphi|^k,\mu\rangle/(k!\beta^k).$ Therefore, for all $k\geq 1,$
$\langle|\varphi|^k,\mu\rangle\leq k!\beta^k.$ It follows that for
all $s\geq 0,$
\begin{eqnarray*}
  \Lambda_\varphi(s) &=& \log \left(1+\sum_{k\geq 1} s^k \langle\varphi^k,\mu\rangle/k!\right) - s\langle\varphi,\mu\rangle \\
   &\leq& \sum_{k\geq 2} s^k \langle\varphi^k,\mu\rangle/k! \\
   &\leq& \sum_{k\geq 2} s^k \langle|\varphi|^k,\mu\rangle/k! \\
   &\leq& \sum_{k\geq 2} (\beta s)^k \\
   &=& \left\{%
\begin{array}{ll}
    (\beta s)^2/(1-\beta s), & \hbox{if } 0\leq \beta s<1 \\
    +\infty, & \hbox{if }\beta s\geq 1 \\
\end{array}%
\right.
\end{eqnarray*}
The last statement holds since $\beta\mapsto\sum_{k\geq 2} (\beta
s)^k $ is an increasing function, for all $s \geq 0.$
\endproof

We are now ready to prove some corollaries of Theorem
\ref{res-BG}.

For any measurable function $f$ in $L_{\rho^*}(\mu),$ let
\begin{eqnarray*}
  \|f\|_\rho^* &:=& \sup\left\{ \IX f\varphi\,d\mu; \varphi: \textrm{measurable, }
\|\varphi\|_\rho\leq 1\right\} \\
   &=& \sup\left\{ \IX f\varphi\,d\mu; \varphi: \textrm{measurable, }
\IX e^{|\varphi|}\,d\mu\leq 2\right\}
\end{eqnarray*}

be the dual norm of $\|\cdot\|_\rho.$

\begin{cor}\label{res-31}
For any probability measure $\nu$ which is absolutely continuous
with respect to $\mu$ and such that $\frac{d\nu}{d\mu}\in
L_{\rho^*}(\mu),$ we have
\[
\left\|\frac{d\nu}{d\mu}-1\right\|_\rho^*\leq 2\sqrt{\Hnm}+\Hnm.
\]
Note that this is the NEI:
$\alpha_1(\|\frac{d\nu}{d\mu}-1\|_\rho^*)\leq \Hnm,$ with
$\alpha_1(t)=(\sqrt{t+1}-1)^2.$
\end{cor}

\proof Here $U$ is the unit ball of
$\mathcal{F}_{\mathrm{exp}}(\mu)$ and thanks to Lemma
\ref{res-Bernstein-Orlicz} applied with $M=1,$ (\ref{eq-35}) holds
as follows: $\Lambda_\varphi(s)\leq\alpha_1^{\mc}(s):=s^2/(1-s).$
Taking the monotone conjugate, we obtain
$\alpha_1(t)=(\sqrt{t+1}-1)^2,$ which is the desired result.
\endproof

The following corollary has already been obtained by F.~Bolley and
C.~Villani in \cite{BV03} with other constants.

\begin{cor}[Weighted CKP inequalities]\label{res-32}
Let $\chi$ be a nonnegative function such that $\IX
e^{a_o\chi}\,d\mu<\infty$ for some $a_o>0.$ Then,
$\|\chi\|_\rho<\infty$ and for any probability measure $\nu$ which
is absolutely continuous with respect to $\mu$ and such that
$\frac{d\nu}{d\mu}\in L_{\rho^*}(\mu),$
$\|\chi\cdot(\nu-\mu)\|_\mathrm{TV}$ is well defined, finite and
we have
\[
\|\chi\cdot(\nu-\mu)\|_\mathrm{TV}\leq \|\chi\|_\rho
\left(2\sqrt{\Hnm}+\Hnm\right)
\]
Note that this is the NEI:
$\alpha(\|\chi\cdot(\nu-\mu)\|_\mathrm{TV})\leq \Hnm,$ with
$\alpha(t)=(\sqrt{t/\|\chi\|_\rho+1}-1)^2.$
\end{cor}

\proof Here $U=\{\chi\psi; \sup|\psi|\leq 1\}.$ As $\chi$ may not
be in $\mathcal{F}_\mathrm{exp}(\mu)$ (if there exists $a_1>0$
such that $\IX e^{a_1\chi}\,d\mu=\infty$), one must be careful. It
happens that
\begin{eqnarray*}
  \|\chi\cdot(\nu-\mu)\|_\mathrm{TV}
  &=& \sup\left\{\IX \chi\psi\, d(\nu-\mu); \psi: \textrm{measurable, }\sup|\psi|\leq 1\right\} \\
   &=& \sup\left\{\IX \varphi\, d(\nu-\mu); \varphi: \textrm{measurable, }|\varphi|\leq \chi,
   \sup|\varphi|<\infty\right\}.
\end{eqnarray*}
To show this, decompose $\nu-\mu$ into its positive and negative
parts, approximate from below $\chi|\psi|
\1_{\mathrm{supp}((\nu-\mu)_+)}$ and $\chi|\psi|
\1_{\mathrm{supp}((\nu-\mu)_-)}$ by  pointwise converging
sequences of bounded functions, and conclude with the dominated
convergence theorem.

Therefore, $U$ can be replaced with
 $U'=\{\varphi; |\varphi|\leq\chi, \sup|\varphi|<\infty\}\subset\mathcal{F}_\mathrm{exp}(\mu).$
As $\sup_{\varphi\in U'}\|\varphi\|_\rho\leq \|\chi\|_\rho,$
thanks to Lemma \ref{res-Bernstein-Orlicz} applied with
$M=\|\chi\|_\rho,$ (\ref{eq-35}) holds as follows:
$\Lambda_\varphi(s)\leq\alpha_M^{\mc}(s):=(Ms)^2/(1-Ms).$ Taking
the monotone conjugate, we obtain
$\alpha_M(t)=(\sqrt{t/M+1}-1)^2,$ which is the desired result.
\endproof

\begin{remark}
It follows from Corollaries \ref{res-31} and \ref{res-32}, that
$$\|\nu-\mu\|_\mathrm{TV}\leq \frac{1}{\log
2}\left(2\sqrt{\Hnm}+\Hnm\right),$$ which of course is worse than
CKP inequality (\ref{eq-Pinsker}) but has the same order of growth
$\sqrt H$ for vanishing entropies.
\end{remark}

Let $d$ be a metric on $\X.$ The associated dual Lipschitz norm of
any signed bounded measure $\xi$ with \emph{zero mass} is defined
by
\[
\|\xi\|_{\mathrm{Lip}}^*=\sup\left\{\IX \varphi\,d\xi; \varphi:
\textrm{measurable, } \|\varphi\|_\mathrm{Lip}\leq 1,
\sup|\varphi|<\infty\right\}
\]
where
$\|\varphi\|_{\mathrm{Lip}}=\sup_{x\not=y}\frac{|\varphi(x)-\varphi(y)|}{d(x,y)}$
is the usual Lipschitz seminorm.

\begin{cor}\label{res-33}
Suppose that there exist $a_o>0$ and $x_o\in\X$ such that $\IX
e^{a_o d(x_o,x)}\,\mu(dx)<\infty.$ Then, $\|d\|_{\rho,\mu^{\otimes
2}}=\inf\{b>0; \int_{\X\times\X} e^{d(x,y)/b}\,\mu(dx)\mu(dy)\leq
2\}<\infty$ and
\[
\left\|\nu-\mu\right\|_\mathrm{Lip}^*\leq \|d\|_{\rho,\mu^{\otimes
2}}\left(2\sqrt{\Hnm}+\Hnm\right),\quad \forall \nu\in\PX.
\]
Note that this is the NEI: $\alpha(\|\nu-\mu\|_\mathrm{Lip}^*)\leq
\Hnm,$ with $\alpha(t)=(\sqrt{t/\|d\|_{\rho,\mu^{\otimes
2}}+1}-1)^2.$
\end{cor}

\proof This is a corollary of Theorem \ref{res-BG-T1}. Here
$U=\{\varphi: \|\varphi\|_\mathrm{Lip}\leq 1,
\sup|\varphi|<\infty\}\subset \mathcal{F}_\mathrm{exp}(\mu).$ Let
us show that
\begin{equation}\label{eq-36}
    \sup_{\varphi\in
U}\|\varphi-\langle\varphi,\mu\rangle\|_\rho\leq
\|d\|_{\rho,\mu^{\otimes 2}}.
\end{equation}
By Jensen's inequality, for any 1-Lipschitz function $\varphi$ and
all $s\geq 0,$
\begin{eqnarray*}
  \exp \left[s\left(\varphi(x)-\IX\varphi(y)\,\mu(dy)\right)\right] &\leq& \IX \exp [s(\varphi(x)-\varphi(y))]\,\mu(dy) \\
   &\leq& \IX \exp [sd(x,y)]\,\mu(dy).
\end{eqnarray*}
Hence, integrating with respect to $\mu(dx),$ one obtains
(\ref{eq-36}).

Thanks to Lemma \ref{res-Bernstein-Orlicz} applied with
$M=\|d\|_{\rho,\mu^{\otimes 2}},$ (\ref{eq-35bis}) holds as
follows: $\Lambda_\varphi(s)\leq\alpha_M^{\mc}(s):=(Ms)^2/(1-Ms).$
Taking the monotone conjugate, we obtain
$\alpha_M(t)=(\sqrt{t/M+1}-1)^2,$ which is the desired result.
\endproof

\section{Tensorization of convex TCIs}\label{sec-tensorization}

In this section only convex TCIs are considered. It is assumed
that the appearing state spaces are Polish and the appearing cost
functions are nonnegative \emph{continuous} and satisfy
(\ref{eq-21}).

\subsection{Statement of the main result}

Let $\mu_1,$ $\mu_2$ be two probability measures on two Polish
spaces $\X_1,$ $\X_2,$ respectively. The cost functions
$c_1(x_1,y_1)$ and
 $c_2(x_2,y_2)$ on $\X_1\times\X_1$ and $\X_2\times\X_2$
give rise to the optimal transportation cost functions
$\mathcal{T}_{c_1}(\mu_1,\nu_1),$ $\nu_1\in\mathcal{P}(\X_1)$ and
$\mathcal{T}_{c_2}(\mu_2,\nu_2),$ $\nu_2\in\mathcal{P}(\X_2).$

On the product space $\X_1\times\X_2,$ we now consider the product
measure $\mu_1\otimes\mu_2$ and the cost function
\[
c_1\oplus
c_2\big((x_1,y_1),(x_2,y_2)\big):=c_1(x_1,y_1)+c_2(x_2,y_2),\quad
x_1,y_1\in\X_1, x_2,y_2\in\X_2
\]
which give rise to the so-called tensorized optimal transportation
cost function
\[
\mathcal{T}_{c_1\oplus c_2}(\mu_1\otimes\mu_2,\nu), \quad\nu\in
\mathcal{P}(\X_1\times\X_2).
\]

Recall that the inf-convolution of two functions $\alpha_1$ and
$\alpha_2$ on $[0,\infty)$ is defined by
\[
\alpha_1 \square \alpha_2(t)=\inf\{\alpha_1(t_1)+ \alpha_2(t_2);
t_1, t_2\geq 0: t_1+t_2=t\},\quad t\geq 0.
\]

\begin{lemma}\label{res-41}
Let $\alpha_1$  and $\alpha_2$ belong to the class $\mathcal{C}.$
Then,
\begin{enumerate}
    \item[(a)]
    $\alpha_1 \square \alpha_2\in\mathcal{C}$ and
    \item[(b)]
    $(\alpha_1 \square \alpha_2)^{\mc}=\alpha_1^{\mc} + \alpha_2^{\mc}$
\end{enumerate}
\end{lemma}
\proof
This simple exercice is left to the reader.
\endproof

The main result of this section is the following theorem.

\begin{thm}[Tensorization]\label{res-tensorization}
Let $c_1$ and $c_2$ be two continuous nonnegative cost functions
which satisfy (\ref{eq-21}).
 Suppose that the convex TCIs
\begin{eqnarray*}
  \alpha_1(\mathcal{T}_{c_1}(\mu_1,\nu_1))
  &\leq& H(\nu_1\mid\mu_1),\ \forall\nu_1\in\mathcal{P}(\X_1) \\
 \alpha_2(\mathcal{T}_{c_2}(\mu_2,\nu_2))
  &\leq& H(\nu_2\mid\mu_2),\ \forall\nu_2\in\mathcal{P}(\X_2)
\end{eqnarray*}
hold with $\alpha_1, \alpha_2\in\mathcal{C}.$ Then, on the product
space $\X_1\times\X_2,$ we have the convex TCI
\begin{equation*}
    \alpha_1 \square \alpha_2\big(\mathcal{T}_{c_1\oplus c_2}(\mu_1\otimes\mu_2,\nu)\big)
    \leq H(\nu\mid\mu_1\otimes\mu_2),\quad\forall \nu\in
\mathcal{P}(\X_1\times\X_2)
\end{equation*}
\end{thm}

Its proof is postponed to Section \ref{sec-dual proof}. We prefer
beginning with a presentation at the next section of an incomplete
derivation of this result which, to our opinion, seems to be more
intuitively appealing.

\subsection{An incomplete direct proof of Theorem \ref{res-tensorization}}

By means of Marton's coupling argument \cite{Mar96}, one can
expect to prove the next Proposition \ref{res-Marton}. We are
interested in transportation costs from $\X_1$ to $\Y_1,$ from
$\X_2$ to $\Y_2$ and from $\X_1\times\X_2$ to $\Y_1\times\Y_2.$

For any probability measure  $\nu$ on the product space
 $\Y=\Y_1\times\Y_2,$ let us write the desintegration of $\nu$
(conditional expectation) as follows:
 $\nu(dy_1dy_2)=\nu_1(dy_1)\nu_2^{y_1}(dy_2).$

 \begin{prop}\label{res-Marton}
For all $\nu\in\mathcal{P}(\Y_1\times\Y_2),$
\begin{equation}\label{eq-41}
     \mathcal{T}_{c_1\oplus c_2}(\mu_1\otimes\mu_2,\nu)
 \leq
 \mathcal{T}_{c_1}(\mu_1,\nu_1)+\int_{\Y_1}\mathcal{T}_{c_2}(\mu_2,\nu_2^{y_1})\,\nu_1(dy_1).
\end{equation}
 \end{prop}

Recall that the relative entropy satisfies for all
$\nu\in\mathcal{P}(\Y_1\times\Y_2),$
\begin{equation}\label{eq-42}
    H(\nu\mid\mu_1\otimes\mu_2)=H(\nu_1\mid\mu_1)+\int_{Y_1}
H(\nu_2^{y_1}\mid\mu_2)\,\nu_1(dy_1)
\end{equation}
which looks like (\ref{eq-41}).

Admitting Proposition \ref{res-Marton} for a while, one can easily
derive Theorem \ref{res-tensorization} as follows. Take
$\Y_1=\X_1$ and $\Y_2=\X_2.$ For all
$\nu\in\mathcal{P}(\X_1\times\X_2),$
\begin{eqnarray*}
 & &  \alpha_1\square\alpha_2 (\mathcal{T}_{c_1\oplus c_2}(\mu_1\otimes\mu_2,\nu))\\
  &\stackrel{(a)}{\leq}&\alpha_1\square\alpha_2\left( \mathcal{T}_{c_1}(\mu_1,\nu_1)+\int_{\Y_1}\mathcal{T}_{c_2}(\mu_2,\nu_2^{y_1})\,\nu_1(dy_1)\right) \\
   &\stackrel{(b)}\leq& \alpha_1(\mathcal{T}_{c_1}(\mu_1,\nu_1))+\alpha_2\left(\int_{\Y_1}\mathcal{T}_{c_2}(\mu_2,\nu_2^{y_1})\,\nu_1(dy_1)\right) \\
    &\stackrel{(c)}\leq& \alpha_1(\mathcal{T}_{c_1}(\mu_1,\nu_1))+\int_{\Y_1}\alpha_2\big(\mathcal{T}_{c_2}(\mu_2,\nu_2^{y_1})\big)\,\nu_1(dy_1) \\
   &\stackrel{(d)}\leq& H(\nu_1\mid\mu_1)+\int_{\Y_1}H(\nu_2^{y_1}\mid\mu_2)\,\nu_1(dy_1) \\
   &=&  H(\nu\mid\mu_1\otimes\mu_2).
\end{eqnarray*}
Inequality (a) holds thanks to Proposition \ref{res-Marton} since
$\alpha_1\square\alpha_2$ is increasing, (b) follows from the very
definition of the inf-convolution, (c) follows from Jensen's
inequality since $\alpha_2$ is convex, (d) follows from the
assumptions $\alpha_1(\mathcal{T}_1(\nu_1))\leq H(\nu_1\mid\mu_1)$
for all $\nu_1$ and $\alpha_2(\mathcal{T}_2(\nu_2))\leq
H(\nu_2\mid\mu_2)$ for all $\nu_2$ (with obvious notations) and
the last equality is (\ref{eq-42}).

To complete the proof of Theorem \ref{res-tensorization}, it
remains to prove Proposition \ref{res-Marton}. This won't be
achieved completely: a difficult measurability statement will only
be conjectured.

\proof[Incomplete proof of Proposition \ref{res-Marton}] One first
faces a nightmare of notations. It might be helpful to introduce
random variables and see
 $\pi\in\mathcal{P}(\X\times\Y)=\mathcal{P}(\X_1\times\X_2\times\Y_1\times\Y_2)$
as the law of $(X_1,X_2,Y_1,Y_2).$ One denotes
 $\pi_1=\mathcal{L}(X_1,Y_1),$
 $\pi_2^{x_1,y_1}\mathcal{L}(X_2,Y_2\mid X_1=x_1, Y_1=y_1),$
 $\pi_{X_2}^{x_1,y_1}=\mathcal{L}(X_2\mid X_1=x_1,Y_1=y_1),$
 $\pi_{Y_2}^{x_1,y_1}=\mathcal{L}(Y_2\mid X_1=x_1,Y_1=y_1)$
 $\pi_X=\mathcal{L}(X_1,X_2),$  $\pi_Y=\mathcal{L}(Y_1,Y_2)$ and so on.

 Let us denote $P(\mu,\nu)$ the set of all $\pi\in\mathcal{P}(\X\times\Y)$ such that
 $\pi_X=\mu$ and $\pi_Y=\nu,$  $P_1(\mu_1,\nu_1)$ the set of all $\eta\in\mathcal{P}(\X_1\times\Y_1)$ such that
 $\eta_{X_1}=\mu_1$ and $\eta_{Y_1}=\nu_1$ and
 $P_2(\mu_2,\nu_2)$ the set of all $\eta\in\mathcal{P}(\X_2\times\Y_2)$ such that
 $\eta_{X_2}=\mu_2$ and $\eta_{Y_2}=\nu_2.$

We only consider couplings $\pi$ such that under the law $\pi$
\begin{itemize}
    \item $\mathcal{L}(X_1,X_2)=\mu,$
    \item $\mathcal{L}(Y_1,Y_2)=\nu,$
    \item $Y_1$ and $X_2$ are independent conditionally on $X_1$ and
    \item $X_1$ and $Y_2$ are independent conditionally on $Y_1.$
\end{itemize}
Optimizing over this collection of couplings leads us to
\[
\mathcal{T}_c(\mu,\nu)\leq \inf_{\pi_1,\pi_2^{\diamond}} \int
c_1\oplus
c_2(x_1,y_1,x_2,y_2)\,\pi_1(dx_1dy_1)\pi_2^{x_1,y_1}(dx_2dy_2)
\]
where the infimum is taken over all $\pi_1\in P_1(\mu_1,\nu_1)$
and all $\pi_2^\diamond=(\pi_2^{x_1,y_1}; x_1\in\X_1, y_1\in\Y_1)$
such that $\pi_2^{x_1,y_1}\in
P_2(\mu_{X_2}^{x_1},\nu_{Y_2}^{y_1})$ for $\pi_1$-almost every
$(x_1,y_1).$ As $\mu$ is a tensor product:
$\mu=\mu_1\otimes\mu_2,$ we have $\mu_{X_2}^{x_1}=\mu_2,$
$\pi_1$-a.e. so that $\pi_2^{x_1,y_1}\in
P_2(\mu_2,\nu_{Y_2}^{y_1})$ for $\pi_1$-almost every $(x_1,y_1).$

Not being careful, one may write
\begin{eqnarray*}
    & & \mathcal{T}_c(\mu,\nu)\\
   &\leq& \inf_{\pi_1,\pi_2^{\diamond}} \int c_1\oplus c_2(x_1,y_1,x_2,y_2)\,\pi_1(dx_1dy_1)\pi_2^{x_1,y_1}(dx_2dy_2) \\
   &=& \inf_{\pi_1} \left[\int_{\X_1\times\Y_1} c_1\,d\pi_1
        + \int_{\X_1\times\Y_1}\left(\inf_{\pi_2^{\diamond}}\int_{\X_2\times\Y_2}c_2(x_2,y_2)\pi_2^{x_1,y_1}(dx_2dy_2)\right)\,\pi_1(dx_1dy_1)\right] \\
   &\stackrel{(a)}{=}& \inf_{\pi_1} \left[\int_{\X_1\times\Y_1} c_1\,d\pi_1
        +
        \int_{\X_1\times\Y_1}\left(\int_{\X_2\times\Y_2}c_2\,d\widehat{\pi}_2^{x_1,y_1}\right)\,\pi_1(dx_1dy_1)\right]\\
  &=& \inf_{\pi_1} \left[\int_{\X_1\times\Y_1} c_1\,d\pi_1
        +
        \int_{\X_1\times\Y_1}\mathcal{T}_{c_2}\big(\mu_2,\nu_{Y_2}^{y_1}\big)\,\pi_1(dx_1dy_1)\right]\\
  &=& \inf_{\pi_1}\left\{\int_{\X_1\times\Y_1} c_1\,d\pi_1\right\}+ \int_{\Y_1}\mathcal{T}_{c_2}(\mu_2,\nu_{Y_2}^{y_1})\,\nu_1(dy_1) \\
   &=& \mathcal{T}_{c_1}(\mu_1,\nu_1)+ \int_{\Y_1}\mathcal{T}_{c_2}(\mu_2,\nu_{Y_2}^{y_1})\,\nu_1(dy_1)
\end{eqnarray*}
which is the desired result.

On the right-hand side of equality (a),
$\widehat{\pi}_2^{x_1,y_1}$ is a minimizer of
$\pi_2^{x_1,y_1}\mapsto\int_{\X_2\times\Y_2}c_2\,d\pi_2^{x_1,y_1}$
subject to the constraint $\pi_2^{x_1,y_1}\in
P_2(\mu_2,\nu_{Y_2}^{y_1}).$ The general theory of optimal
transportation insures that such a minimizer exists for each
$(x_1,y_1).$ And it might seem that the work is done.

But this is not true since one still has to prove that there
exists a \emph{measurable} mapping $(x_1,y_1)\mapsto
\widehat{\pi}_2^{x_1,y_1}.$ We now face a difficult problem that
may possibly be solved by means of a measurable selection theorem,
taking advantage of the pleasant property of tightness of any
probability measure on a Polish space.

We withdraw this promising direct approach.
\endproof

\subsection{A complete indirect proof of Theorem \ref{res-tensorization}}\label{sec-dual proof}

It is based upon an indirect dual approach, making use of the
characterization of Corollary \ref{res-34} and follows the line of
proof of (\cite{Led01}, Proposition 1.19).

\proof[Proof of Theorem \ref{res-tensorization}]

 Recall that, provided that $c$ is continuous nonnegative and satisfy (\ref{eq-21}), $Q^c\varphi (x)=\inf_{y\in\X}\{\varphi(y)+c(y,x)\}$
 is in $B(\X)$ whenever
 $\varphi\in B(\X).$  We denote $Q_1=Q^{c_1},$
 $Q_2=Q^{c_2}$ and $Q=Q^{c_1\oplus c_2}.$

 By Corollary \ref{res-34}, the convex TCIs ``$\alpha_1(\mathcal{T}_1)\leq
 H_1$" and ``$\alpha_2(\mathcal{T}_2)\leq
 H_2$" which are suppposed to hold are equivalent to
 \begin{eqnarray}
   \int_{\X_1} e^{sQ_1\theta_1}\,d\mu_1&=& \exp(\alpha_1^{\mc}(s)+s\langle\theta_1,\mu_1\rangle),
   \quad \forall s\geq 0, \forall \theta_1\in B(\X_1) \label{eq-43} \\
   \int_{\X_2} e^{sQ_2\theta_2}\,d\mu_2&=& \exp(\alpha_2^{\mc}(s)+s\langle\theta_2,\mu_2\rangle),
   \quad \forall s\geq 0, \forall \theta_2\in B(\X_2)
    \label{eq-44}
 \end{eqnarray}

As by Lemma \ref{res-41} $(\alpha_1 \square
\alpha_2)^{\mc}=\alpha_1^{\mc} + \alpha_2^{\mc},$ thanks to
Corollary  \ref{res-34} again, all we have to prove is
\begin{equation}\label{eq-45}
     \int_{\X_1\times\X_2} e^{sQ\varphi}\,d(\mu_1\otimes \mu_2)
     = \exp(\alpha_1^{\mc}+\alpha_2^{\mc}(s)+s\langle\varphi,\mu_1\otimes \mu_2\rangle),
   \quad \forall s\geq 0, \forall \varphi\in C_b(\X_1\times\X_2)
\end{equation}

Let us take $\varphi\in C_b(\X_1\times\X_2).$ For all
$(x_1,x_2)\in\X_1\times\X_2,$
\begin{eqnarray*}
  Q\varphi(x_1,x_2) &=& \inf_{y_1\in\X_1, y_2\in Y_2} \{\varphi(y_1,y_2)+c_1(y_1,x_1)+c_2(y_2,x_2)\} \\
   &=& \inf_{y_1\in\X_1}\left\{ \inf_{y_2\in Y_2}\{\varphi(y_1,y_2)+c_2(y_2,x_2)\}+c_1(y_1,x_1)\right\} \\
   &=& \inf_{y_1\in\X_1}\{ \theta_{x_2}(y_1)+c_1(y_1,x_1)\}\\
   &=& Q_1\theta_{x_2}(x_1)
\end{eqnarray*}
where
\begin{equation}\label{eq-46}
    \theta_{x_2}(y_1)=Q_2\varphi_{y_1}(x_2)=\inf_{y_2\in
Y_2}\{\varphi(y_1,y_2)+c_2(y_2,x_2)\}
\end{equation}
with $\varphi_{y_1}(y_2):=\varphi(y_1,y_2).$ Hence, for all $s\geq
0,$
\begin{eqnarray*}
    \int_{\X_1\times\X_2} e^{sQ\varphi}\,d(\mu_1\otimes \mu_2)
    &\stackrel{(a)}{=}& \int_{\X_2}\left(\int_{\X_1}e^{sQ_1\theta_{x_2}(x_1)}\,\mu_1(dx_1)\right)\,\mu_2(dx_2) \\
   &\stackrel{(b)}{\leq}& \int_{\X_2}e^{\alpha_1^{\mc}(s)+s\langle \theta_{x_2},\mu_1\rangle}\,\mu_2(dx_2) \\
   &\stackrel{(c)}{=}& e^{\alpha_1^{\mc}(s)} \int_{\X_2}
   \exp\left(s\int_{\X_1}Q_2\varphi_{y_1}(x_2)\,\mu_1(dy_1)\right)\,\mu_2(dx_2)
\end{eqnarray*}
Equality (a) is justified since $\varphi$ being bounded,
$(x_1,x_2)\mapsto Q\varphi(x_1,x_2)=Q_1\theta_{x_2}(x_1)$ is
jointly measurable.
\\
Let us now prove the inequality (b). As $\varphi$ and $c$ are
continuous, $(x_2,y_1)\mapsto\theta_{x_2}(y_1)$ is jointly upper
semicontinuous as the infimum of a collection of continuous
functions. Since $ \theta_{x_2}(y_1)=Q_2\varphi_{y_1}(x_2)$ by
(\ref{eq-46}), we have $\sup_{x_1,x_2}|\theta_{x_2}(y_1)|\leq
\sup_{y_1}\sup |\varphi_{y_1}|= \sup |\varphi|<\infty.$ Therefore,
$(x_2,y_1)\mapsto\theta_{x_2}(y_1)$ is an upper semicontinuous
bounded function. Consequently, one is allowed to invoke
(\ref{eq-43}) to obtain
$\int_{\X_1}e^{sQ_1\theta_{x_2}(x_1)}\,\mu_1(dx_1)\leq
e^{\alpha_1^{\mc}(s)+s\langle \theta_{x_2},\mu_1\rangle}$ for all
$x_2.$ Also note that $x_2\mapsto \langle
\theta_{x_2},\mu_1\rangle$ is measurable since $(x_2,y_1)\mapsto
\theta_{x_2}(y_1)$ is jointly measurable and bounded.
\\
The last equality (c) is simply (\ref{eq-46}).

\begin{remark}
 If $c_2$ is only assumed to be \lsc, the joint
measurability of $(x_2,y_1)\mapsto\theta_{x_2}(y_1)$ which has
been used to prove inequality (b) is far from being clear. This is
the reason why the cost functions are supposed to be continuous.
\end{remark}

But for all $x_2,$
\begin{eqnarray*}
  \int_{\X_1}Q_2\varphi_{y_1}(x_2)\,\mu_1(dy_1)
  &=&  \int_{\X_1}\inf_{y_2\in Y_2}\{\varphi(y_1,y_2)+c_2(y_2,x_2)\}\,\mu_1(dy_1)\\
   &\leq& \inf_{y_2\in Y_2}\left\{\int_{\X_1}\varphi(y_1,y_2)\,\mu_1(dy_1)+c_2(y_2,x_2)\right\} \\
   &=& Q_2\overline{\varphi}(x_2)
\end{eqnarray*}
where
$y_2\mapsto\overline{\varphi}(y_2)=\int_{\X_1}\varphi(y_1,y_2)\,\mu_1(dy_1)$
is a continuous bounded function. Gathering our partial results
leads us, for all $s\geq 0,$ to the inequality (a) below
\begin{eqnarray*}
   \int_{\X_1\times\X_2} e^{sQ\varphi}\,d(\mu_1\otimes \mu_2)
   &\stackrel{(a)}{\leq}& e^{\alpha_1^{\mc}(s)} \int_{\X_2} e^{sQ_2\overline{\varphi}}\,d\mu_2 \\
   &\stackrel{(b)}{\leq}& e^{\alpha_1^{\mc}(s)} e^{\alpha_2^{\mc}(s)+s\langle\overline{\varphi},\mu_2\rangle} \\
   &=& e^{\alpha_1^{\mc}(s)+\alpha_2^{\mc}(s)+s\langle\varphi,\mu_1\otimes\mu_2\rangle}
\end{eqnarray*}
Inequality (b) is a consequence of (\ref{eq-44}). This is
(\ref{eq-45}) and concludes the proof of the theorem.
\endproof

\subsection{Product of $n$ spaces}
The extension of Theorem \ref{res-tensorization} to the product of
$n$ spaces is as follows. Let $\X_1,\dots,\X_n$ be $n$ Polish
spaces and $\mu_1,\dots,\mu_n$ be probability measures on each of
these spaces. On each space $\X_i$ let $c_i$ be a cost function.
The cost function on the product space
$\X_1\times\cdots\times\X_n$ is
\[
c_1\oplus\cdots\oplus
c_n\big((x_1,\dots,x_n),(y_1,\dots,y_n)\big)=c_1(x_1,y_1)+\cdots
+c_n(x_n,y_n)
\]

\begin{cor}\label{res-n-tensorization}
Let us assume that the cost functions $c_i$ are nonnegative
continuous and satisfy (\ref{eq-21}). Suppose that the convex
transportation cost inequalities
\begin{equation*}
  \alpha_i(\mathcal{T}_{c_i}(\mu_i,\nu_i))
  \leq H(\nu_i\mid\mu_i),\ \forall\nu_i\in\mathcal{P}(\X_i), \quad
  i=1,\dots,n
\end{equation*}
hold with $\alpha_1,\dots, \alpha_n\in\mathcal{C}.$ Then, on the
product space $\X_1\times\cdots\times\X_n,$ we have the convex
transportation cost inequality
\begin{equation*}
  \alpha_1\square\cdots\square\alpha_n  \big(\mathcal{T}_{c_1\oplus\cdots\oplus c_n}(\mu_1\otimes\cdots\otimes\mu_n,\nu)\big)
    \leq H(\nu\mid\mu_1\otimes\cdots\otimes\mu_n),\ \forall \nu\in
\mathcal{P}(\X_1\times\cdots\times\X_n)
\end{equation*}
where
\[
\alpha_1\square\cdots\square\alpha_n(t)=\inf\{\alpha_1(t_1)+\cdots+\alpha_n(t_n);
t_1,\dots,t_n\geq 0: t_1+\cdots+t_n=t\},\  t\geq 0
\]
is the inf-convolution of $\alpha_1,\dots,\alpha_n.$
\end{cor}
\proof It is a direct consequence of Theorem
\ref{res-tensorization} which is proved by induction, noting that
$\alpha_1\square\cdots\square\alpha_n=(\alpha_1\square\cdots\square\alpha_{n-1})\square\alpha_n$
for all $n.$
\endproof

In the special situation where the $n$ TCIs are copies of a unique
TCI on a Polish space $\X$ we have the following important result.

\begin{thm}\label{res-n tenso bis}
Let us assume that the cost function $c$ is nonnegative continuous
and satisfy (\ref{eq-21}). Suppose that the convex transportation
cost inequality
\begin{equation*}
  \alpha(\mathcal{T}_{c}(\mu,\nu))
  \leq H(\nu\mid\mu),\ \forall\nu\in\mathcal{P}(\X)
\end{equation*}
holds with $\alpha \in\mathcal{C}.$ Then, on the product space
$\X^n,$ we have the following convex transportation cost
inequality
\begin{equation*}
n\alpha \left(\frac{\mathcal{T}_{c^{\oplus n}}(\mu^{\otimes
n},\zeta)}{n}\right)
    \leq H(\zeta\mid\mu^{\otimes n}),\ \forall \zeta\in
\mathcal{P}(\X^n)
\end{equation*}
where $c^{\oplus
n}\big((x_1,\dots,x_n),(y_1,\dots,y_n)\big)=c(x_1,y_1)+\cdots
+c(x_n,y_n).$
\end{thm}
\proof This is a direct application of Corollary
\ref{res-n-tensorization}, noting that $\alpha^{\square
n}(t)=n\alpha(t/n).$
\endproof

\par\medskip
\textbf{About dimension-free tensorized convex TCIs.}\ Let us say
that a convex transportation cost inequality
\begin{equation}\label{eq-47}
\alpha\left(\mathcal{T}_{c}(\mu,\nu)\right)\leq H(\nu\mid\mu),\quad \forall \nu\in \mathcal{P}(\X)
\end{equation}
has the dimension-free tensorization property, if the inequality
$$\alpha\left(\mathcal{T}_{c^{\oplus n}}(\mu^{\otimes n},\zeta)\right)\leq H(\zeta\mid\mu^{\otimes n}),\quad \forall \zeta\in \mathcal{P}(\X^n)$$
holds for all $n\in \mathbb{N}^*$.\medskip\\
Clearly, according to Theorem \ref{res-n tenso bis}, if $\alpha\in \mathcal{C}$ is of the form $\alpha(t)=at$
with $a\geq0$, then (\ref{eq-47}) has the dimension-free tensorization property.

\begin{remark}
Thanks to the same theorem, a seemingly weaker sufficient
condition on $\alpha$ for (\ref{eq-47}) to be dimension-free is
$\alpha(t)\leq \inf_{n\geq 1}n\alpha(t/n),$ $t\geq 0.$ As $\alpha$
is in $\mathcal{C},$ $\alpha(t)/t$ is an increasing function so
that $\alpha'(0):=\lim_{t\downarrow 0}\alpha(t)/t$ exists. It
follows that $\lim_{n\rightarrow\infty} n\alpha(t/n)=\alpha'(0)t$
for all $t\geq 0.$ Therefore, the condition $\alpha(t)\leq
\inf_{n\geq 1}n\alpha(t/n),$ $t\geq 0$ is equivalent to
$\alpha(t)\leq \alpha'(0)t,$ $t\geq 0$. But since $\alpha$ is
convex, the converse inequality also holds, that is $\alpha(t)\geq
\alpha'(0)t,$ $t\geq 0$. Consequently $\alpha$ is of the form
$\alpha(t)=at$ with $a\geq 0$.
\end{remark}

Dimension free tensorization is a phenomenon that can only happen when dealing with \emph{non-metric} cost functions. Indeed, we show in the following proposition, that convex $T_1$-inequalities having this property are all trivial.
\begin{prop}
Let $(\X,d)$ be a Polish space and $\mu\in \PX$. The convex transportation cost inequality
\begin{equation}\label{eq-48}
\alpha\left(\mathcal{T}_{d}(\mu,\nu)\right)\leq H(\nu\mid\mu),\quad \forall \nu\in \mathcal{P}(\X),
\end{equation}
with $\alpha\in \mathcal{C}$ has the dimension free tensorization property if, and only if $\alpha=0$ or $\mu$ is a Dirac mass.
\end{prop}
\proof If $\alpha=0$, it is clear that (\ref{eq-48}) has the dimension free tensorization property. If $\mu$ is a Dirac mass, it is easy to see that (\ref{eq-48}) holds for every $\alpha \in \mathcal{C}$. Noting that a tensor product of Dirac measures is again a Dirac measure, the dimension-free tensorization property is established in this special case.\\
Now, suppose that (\ref{eq-48}) has the dimension-free tensorization property, with $\alpha\neq 0$ and let us prove that $\mu$ is a Dirac mass.
According to Theorem \ref{res-BG-T1}, the following inequality
$$\log\int_{\X^n}e^{s(\varphi(x_1)+\cdots+\varphi(x_n)-n\langle\varphi,\mu\rangle)}\,\mu^{\otimes n}(dx_1\ldots dx_n)\leq \alpha^\circledast(s),\quad\forall s\geq0 $$
holds for all bounded $1$-Lipschitz $\varphi$ and all $n\geq 1$. As a consequence, denoting by $\Lambda_\varphi$ the Log-Laplace of $\varphi(X)-\langle\varphi,\mu\rangle$, $X$ of law $\mu$, one has
$\Lambda_\varphi\leq \frac{1}{n}\alpha^\circledast,$ for all $n\geq 1$, and so $\Lambda_\varphi\leq 0$ on $\dom \alpha^\circledast$ (the effective domain of $\alpha^\circledast$). But by Jensen inequality, one obtains immediately $\Lambda_\varphi \geq 0$. Thus $\Lambda_\varphi\equiv0$ on $\dom \alpha^\circledast$. As $\alpha\neq 0$, $[0,a[ \subset\dom \alpha^\circledast$, for some $a>0$. Considering $-\varphi$ instead of $\varphi$ in the above reasoning yields that $\Lambda_\varphi\equiv 0$ on $]-a,a[$. This easily implies that $\mu_\varphi$ (the image of $\mu$ under the application $\varphi$) is a Dirac mass. Now, let us take a point $x_0$ in the support of $\mu$ and consider the bounded $1$-Lipschitz function $\varphi_0(x)=d(x,x_0)\wedge1$, $x\in \X$. As $x_0$ is in the support of $\mu$, $\mu_{\varphi_0}([0,\varepsilon[)=\mu(\varphi_0<\varepsilon)>0$ for all $\varepsilon>0$. As $\mu_{\varphi_0}$ is a Dirac mass, one thus has $\mu(\varphi_0<\varepsilon)=1$ for all $\varepsilon>0$. This easily implies that $\mu=\delta_{x_0}$.
\endproof

\section{Integral criteria}\label{sec-integral criteria}

Our aim in this section is to give integral criteria for a convex
$\mathcal{T}$-inequality to hold.

Let us first note that when two $\mathcal{T}$-inequalities
$\alpha_0(\mathcal{T}(\nu))\leq H(\nu\mid\mu), $ $\forall\nu\in\N$
and $\alpha_1(\mathcal{T}(\nu))\leq H(\nu\mid\mu), $
$\forall\nu\in\N$ hold, then we have the resulting new inequality
$\alpha(\mathcal{T}(\nu))\leq H(\nu\mid\mu), $ $\forall\nu\in\N$
with
\begin{equation}\label{eq-54}
    \alpha=\max(\alpha_0,\alpha_1).
\end{equation}

This allows us to separate our investigation into two parts:
obtaining $\alpha_0$ and $\alpha_1$ which control respectively the
small (neighbourhood of $t=0$) and large values of $t$ (the other
ones). Let us go on with some vocabulary.

\subsection{Transportation functions and deviation functions}

 We introduce the following definitions. Recall
that $\mathcal{T}$ is defined at (\ref{def-T}).

\begin{definition}[Transportation function]\label{def-trans function}
A left continuous increasing function
$\alpha:[0,\infty)\rightarrow [0,\infty]$ is called a
\emph{transportation function} for $\mathcal{T}$ in $\N$ if
\begin{equation*}
    \alpha(\mathcal{T}(\nu))\leq H(\nu\mid \mu),\quad \forall\nu\in\N.
\end{equation*}
\end{definition}

This means that the $\mathcal{T}$-inequality (\ref{eq-TIneq1})
holds with $\alpha.$

\begin{definition}[Deviation function]\label{def-dev function}
A left continuous increasing function
$\alpha:[0,\infty)\rightarrow [0,\infty]$ is called a
\emph{deviation function} for $\mathcal{T}$ if
\begin{equation*}
    \logsup \P(\mathcal{T}(L_n)\geq t)\leq -\alpha(t),\quad\forall
    t\geq 0.
\end{equation*}
\end{definition}

These functions will be shortly called later transportation and
deviation functions, without any reference to $\mathcal{T}$ and
$\N.$

\begin{remark}\label{rem-51}
For $\mathcal{T}(L_n)$ to be measurable, it is assumed that $\Phi$
is a set of couples of \emph{continuous} functions. Indeed,
$$
\left\{\mathbf{x}\in\X^n;
\mathcal{T}\left(\frac{1}{n}\sum_{i=1}^n\delta_{x_i}\right)\leq
t\right\}=\bigcap_{\phi\in\Phi}\left\{\mathbf{x}\in\X^n;\frac{1}{n}\sum_{i=1}^n
\varphi(x_i)+\langle\psi,\mu\rangle\leq t\right\}
$$
is a closed set.
\end{remark}

Note that an increasing function is left continuous if and only if
it is \lsc. Clearly, the best transportation function is the left
continuous version of the increasing function
\[
t\mapsto \inf\{H(\nu\mid\mu); \nu\in\N, \mathcal{T}(\nu)\geq t\},
t\geq 0.
\]
Similarly, the best deviation function is the left continuous
version of the increasing function
\[
 t\mapsto-\logsup\P(\mathcal{T}(L_n)\geq t)\in[0,\infty],
t\geq 0.
\]

\begin{prop}\label{res-51}
Under the assumptions of Theorem \ref{res-BG}, any deviation
function $\alpha$ in the class $\mathcal{C}$ is a transportation
function.
\end{prop}

\proof Let $\alpha\in\mathcal{C}$ be a deviation function. Since
$\mathcal{T}(L_n)\geq T_n^\phi$ for all $\phi\in\Phi,$ we clearly
have $\P(\mathcal{T}(L_n)\geq t)\geq \P(T_n^\phi\geq t)$ for all
$t\geq 0$ and $n.$ Therefore, for all $\phi, n$ and $t,$ $ \logsup
\P(T_n^\phi\geq t)\leq\logsup \P(\mathcal{T}(L_n)\geq t)\leq
-\alpha(t).$ This implies the statement (d) of Theorem
\ref{res-BG}, which in turn is equivalent to the statement (a) of
Theorem \ref{res-BG}, which is the desired result.
\endproof

\subsection{Controlling the large values of $t$}

In this subsection, it is assumed that the deviation and
transportation functions are in $\mathcal{C}.$

\begin{prop}\label{res-53}
The first statement is concerned with convex TCIs and the second
one with convex $\mathcal{T}$-inequalities.
\begin{enumerate}
\item[(a)] If $\beta\in\mathcal{C}$ satisfies $\IX\exp [\beta(\IX
c(x,y)\,\mu(dy))]\,\mu(dx)\leq A<\infty$ then
\[
\alpha(t)=\max(0,\beta(t)-\log A)),\quad t\geq 0
\]
is a transportation function.

\item[(b)] Let us suppose that $\alpha$ is a transportation
function, then for all $\pf\in\Phi$
\[
 \IX \exp \left[\delta
\alpha\left(\varphi(x)+\langle\psi,\mu\rangle\right)\right]\,\mu(dx)
\leq \frac{1+\delta}{1-\delta}<\infty,\quad \forall 0\leq\delta<1.
\]
\end{enumerate}
\end{prop}

\textit{Remarks}.
 \begin{itemize}
    \item In (a), because of Jensen's inequality, one can take $A\geq \int_{\X^2}\exp \beta(c(x,y))\,\mu(dx)\mu(dy)$
    \item About (a), if $c=d\leq D<\infty$ is a \lsc\ bounded metric,
    one recovers that $\alpha(t)=\left\{%
\begin{array}{ll}
    0, & \hbox{if }t\leq D \\
    +\infty, & \hbox{if }t>D \\
\end{array}%
\right.    $ is a transportation function, which is obvious.
    \item About (b) in the case of a TCI, let us note that $\sup_{\pf\in\Phi_c}(\varphi(x)+\langle\psi,\mu\rangle)\leq
    \IX \sup_\phi(\varphi(x)+\psi(y))\,\mu(dy)\leq\IX
    c(x,y)\,\mu(dy)$ for all $x.$
    It follows that
    \\ $\IX \exp \left[\delta
\alpha\left((\varphi(x)+\langle\psi,\mu\rangle)\right)\right]\,\mu(dx)\leq
 \IX \exp \left[\delta
\alpha\left(\IX c(x,y)\,\mu(dy)\right)\right]\,\mu(dx)$ for all
$\pf\in\Phi.$
    It would be pleasant to obtain the finiteness of
    an integral in terms of $c.$ In the case where $c(x,y)=d(x,y)^p,$ this will be performed below at
    Corollary \ref{res-54}.
 \end{itemize}

\proof Let us prove (a). As the product measure $\mu(dx) L_n(dy)$
has the right marginal measures, we get: $\mathcal{T}_c(\mu,
L_n):=T_n\leq \int_{\X^2}c(x,y)\mu(dx)L_n(dy)=\langle
c_\mu,L_n\rangle$ with $c_\mu(y):=\IX c(x,y)\,\mu(dx).$ It follows
that for all $t\geq 0,$
\begin{eqnarray*}
  \P(T_n\geq t) &\leq& \P(\langle c_\mu,L_n\rangle\geq t) \\
   &\stackrel{(a)}{=}& \P(\beta(\langle c_\mu,L_n\rangle)\geq \beta(t)) \\
   &\stackrel{(b)}{\leq}& \P(\langle \beta\circ c_\mu,L_n\rangle\geq \beta(t)) \\
   &\stackrel{(c)}{=}& \P(e^{\sum_{i=1}^n \beta\circ c_\mu(X_i)}\geq e^{n\beta(t)})  \\
  &\stackrel{(d)}{\leq}& e^{-n\beta(t)}\E e^{\sum_{i=1}^n \beta\circ c_\mu(X_i)} \\
   &\stackrel{(e)}{=}& \left[e^{-\beta(t)}\E e^{\beta\circ c_\mu(X)}\right]^n
\end{eqnarray*}
where equality (a) follows from the monotony of $\beta,$ (b) from
the convexity of $\beta$ and Jensen's inequality, (c) from the
monotony of the exponential, (d) from Markov's inequality and (e)
from the fact that $(X_i)$ is an iid sequence. Finally,
\[
\logsup \P(T_n\geq t)\leq -\beta(t)+\log\IX e^{\beta\circ
c_\mu}\,d\mu, \ \forall t\geq 0
\]
which with Proposition \ref{res-51} leads to the desired result.
\par\medskip
Let us prove (b). As $\alpha\in\mathcal{C}$ is a transportation
function, by Theorem \ref{res-BG} (keeping the notations of
Theorem \ref{res-BG}) we have for all
\[
\alpha(t)\leq \Lambda_\phi^*(t),\ \forall\phi\in\Phi, \forall
t\geq 0.
\]
By Lemma \ref{res-52} below, as $\Lambda_\phi^*$ is the Cram\'er
transform of $\varphi(X)+\langle\psi,\mu\rangle$ we get
\begin{equation*}
    \E \exp\left[\delta
    \Lambda_\phi^*(\varphi(X)+\langle\psi,\mu\rangle)\right]\leq
    \frac{1+\delta}{1-\delta},
    \forall 0\leq\delta<1, \forall \phi
\end{equation*}
Extending $\alpha$ with $\alpha(t)=0$ for all $t\leq 0,$ we obtain
$\alpha\leq\Lambda_\phi^*$  for all $\phi,$ $\alpha\leq
\Lambda_\phi^*.$ Consequently we obtain
\begin{equation*}
    \IX \exp\left[\delta
    \alpha(\varphi(x)+\langle\psi,\mu\rangle)\right]\,\mu(dx)\leq
    \frac{1+\delta}{1-\delta},
    \forall 0\leq\delta<1, \forall \phi
\end{equation*}
As $e^{\delta\alpha}$ is increasing, the desired result follows by
monotone convergence.
\endproof

During the above proof, the following lemma has been used.

\begin{lemma}\label{res-52}
Let $Z$ be a  real random variable such that $\E e^{\lambda_o
|Z|}<\infty$ for some $\lambda_o>0.$ Let $h$ be its Cram\'er
transform. Then for all $0\leq \delta<1,$ $\E \exp[\delta
h(Z)]\leq (1+\delta)/(1-\delta).$
\end{lemma}

\proof This result with the upper bound $2/(1-\delta)$ instead of
$(1+\delta)/(1-\delta)$ can be found in (\cite{DZ}, Lemma 5.1.14).
For a proof of the improvement with $(1+\delta)/(1-\delta)$ see
\cite{GozPhD}.
\endproof

\begin{cor}\label{res-55}
In this statement $d$ is a \lsc\ semimetric and $c$ is a \lsc\
cost function such that $c(x,x)=0$ for all $x\in\X.$
\begin{enumerate}
    \item[(a)] Suppose that there exists a nonnegative measurable function
    $\chi$ such that
    $$
    c\leq \chi \oplus\chi.
    $$
    Let $\gamma\in\mathcal{C}$ be such that $\IX \exp[\gamma\circ\chi(x)]\,\mu(dx)\leq
    B<\infty,$ then for any $x_o\in\X$
    \[
t\mapsto 2\max(0,2\gamma(t/4)-\gamma\circ\chi(x_o)-\log B), \quad
t\geq 0
    \]
    is a transportation function for $c.$
    \item[(b)]
    Suppose that there exists $\theta\in\mathcal{C}$ such that
    $$
    \theta(d)\leq c.
    $$
    If $\alpha\in\mathcal{C}$ is a transportation function for $c,$ then
    \[
\IX \exp[u\ \alpha\circ\theta(d(x_o,x)/2)]\,\mu(dx)<\infty
    \]
    for all $x_o\in\X$ and all $0\leq u<2.$
\end{enumerate}
\end{cor}

\proof We begin with the case where $c=d,$ $\chi(x)=d(x_o,x)$ and
$\theta(d)=d.$
\par\medskip
 \textit{The case $c=d$.}\ To prove (a) with $\chi(x)=d(x_o,x),$
we apply  statement (a) of Proposition \ref{res-53}. Let $\beta$
be in the class $\mathcal{C}.$ We have for all $x_o\in\X$
\begin{eqnarray*}
  \IX\exp \Big[\beta\Big(\IX d(x,y)\,\mu(dy)\Big)\Big]\,\mu(dx)
  &\leq& \int_{\X^2} \exp[\beta(d(x,y))]\,\mu(dx)\mu(dy) \\
   &\leq& \int_{\X^2} \exp\Big[\beta\Big(\frac{2d(x_o,x)+2d(x_o,y)}{2}\Big)\Big]\,\mu(dx)\mu(dy) \\
   &\leq& \int_{\X^2} \exp[\beta(2d(x_o,x))/2+\beta(2d(x_o,y))/2]\,\mu(dx)\mu(dy) \\
   &=& \left(\IX \exp\left[\frac{\beta(2d(x_o,x))}{2}\right]\,\mu(dx)\right)^2\\
   &:=&A
\end{eqnarray*}
Taking, $\gamma(t)=\beta(2t)/2,$ one gets $A=B^2$ and
\begin{equation}\label{eq-51}
    t\mapsto\max(0,\beta(t)-\log A)=\max(0,2\gamma(t/2)-2\log B)
\end{equation}
is a transportation function for $c=d.$

\par\medskip
Now, let us prove (b). Thanks to Kantorovich-Rubinstein
equality (\ref{eq-KR}) one can take  $\Phi=\{(-\varphi,\varphi);
\|\varphi\|_{\textrm{Lip}}\leq 1, \varphi \textrm{ bounded} \}.$
Because of Proposition \ref{res-53}-(b), we have for all bounded
$\varphi$ with $\|\varphi\|_{\textrm{Lip}}\leq 1:$
\begin{equation*}
    \IX \exp[\delta\alpha
    (\varphi(x)-\langle\varphi,\mu\rangle)]\,\mu(dx)\leq
    (1+\delta)/(1-\delta), \forall 0\leq\delta<1.
\end{equation*}
 The function $\varphi(x)=d(x_o,x)$ is 1-Lipschitz but it is not
 bounded in general. Let us introduce an approximation procedure.
 For all $k\geq 0,$ with $m:=\IX d(x_o,y)\,\mu(dy),$ we have
 \begin{eqnarray*}
   \IX \exp[\delta\alpha((d(x_o,x)\wedge k) -m]\,\mu(dx)
    &\leq& \IX \exp\left[\delta\alpha
    ((d(x_o,x)\wedge k) -\IX [d(x_o,y)\wedge k]\,\mu(dy)\right]\,\mu(dx) \\
    &\leq&  (1+\delta)/(1-\delta).
 \end{eqnarray*}
By monotone convergence, one concludes that for all
$0\leq\delta<1,$
\[
\IX \exp[\delta\alpha(d(x_o,x) -m]\,\mu(dx)
    \leq  (1+\delta)/(1-\delta).
\]
As
\begin{eqnarray*}
  2\delta\alpha(d(x_o,x)/2) &=& 2\delta\alpha\left(\frac{d(x_o,x)-m}{2}+\frac{m}{2}\right) \\
   &\leq& \delta [\alpha(d(x_o,x)-m)+\alpha(m)],
\end{eqnarray*}
one sees that
\[
\IX \exp[2\delta\alpha(d(x_o,x)/2)]\,\mu(dx)\leq
e^{\delta\alpha(m)}\IX \exp[\delta\alpha(d(x_o,x) -m]\,\mu(dx)
    \leq  e^{\delta\alpha(m)}(1+\delta)/(1-\delta)
\]
which leads to
\begin{equation}\label{eq-53}
    \IX \exp[2\delta \alpha(d(x_o,x)/2)]\,\mu(dx)<\infty
\end{equation}

\par\medskip
 \textit{The general case.}\
Let us prove (a). It is clear that $c(x,y)\leq d_\chi(x,y)$ where
$d_\chi$ is the semimetric defined by
\begin{equation}\label{eq-dchi}
    d_\chi(x,y)=\1_{x\not =y}(\chi(x)+\chi(y)).
\end{equation}

\begin{remark}
If $\chi$ admits two or more zeros, $d_\chi$ is a
semimetric. Otherwise it is a metric. In the often studied case
where $c=d^p$ with $d$ a metric and $p\geq 1,$ one takes
$\chi(x)=2^{p-1}d(x_o,x)^p$ (see the proof of Corollary
\ref{res-54} below) and $d_\chi$ is a metric.
\end{remark}

Of course, for all $\nu\in\N=\mathcal{P}_\chi=\{\nu\in\PX; \IX
\chi(x)\,\nu(dx)<\infty\},$ we have
\[
\mathcal{T}_{c}(\nu)\leq \mathcal{T}_{d_\chi}(\nu).
\]
Therefore, any transportation function for $d_\chi$ is a
transportation function for $c.$ This easy but powerful trick is
borrowed from the monograph by  C.~Villani (\cite{Vill},
Proposition 7.10).
\\
It has been proved at (\ref{eq-51}) that if $\IX
\exp[\beta(d_\chi(x_o,x)]\,\mu(dx)\leq C<\infty$ for some function
$\beta\in\mathcal{C},$ then $\max(0,2\beta(t/2)-2\log C)$ is a
transportation function for $d_\chi.$

Taking $\beta(t)=2\gamma(t/2),$ with convexity we have
\begin{equation}\label{eq-52}
    \beta(d_\chi(x_o,x))\leq \gamma\circ\chi (x_o)+\gamma\circ\chi
    (x)
\end{equation}
so that $\IX \exp[\beta(d_\chi(x_o,x)]\,\mu(dx)\leq
e^{\gamma\circ\chi (x_o)} B=C.$ This leads us to
$\max(0,2\beta(t/2)-2\log
C)=2\max(0,2\gamma(t/4)-\gamma\circ\chi(x_o)-\log B)$ which is the
desired result.

\par\medskip
Let us prove (b). Because of Jensen's inequality, it is easy to
show that $\theta(\mathcal{T}_d)\leq \mathcal{T}_{c}.$ As $\alpha$
is a transportation function for $c,$ it follows that
$\alpha\circ\theta$ is a transportation function for
$\mathcal{T}_d.$ Applying the already proved result (\ref{eq-53})
with  $\alpha\circ\theta$ instead of $\alpha$  completes the proof
of the corollary.
\endproof

Now, we consider an important special  case of convex TCI.

\begin{cor}[$c=d^p$]\label{res-54}
In this statement $c=d^p$ where $d$ is a \lsc\ metric and $p\geq
1.$
\begin{enumerate}
    \item[(a)] Let $\gamma\in\mathcal{C}$ be such that $\IX \exp[\gamma(d^p(x_o,y))]\,\mu(dy)\leq
    B<\infty$ for some $x_o\in\X,$ then
    \[
t\mapsto\max(0,2\gamma(2^{-p}t)-2\log B), \quad t\geq 0
    \]
    is a transportation function.
    \item[(b)]
    If $\alpha\in\mathcal{C}$ is a transportation function, then
    \[
\IX \exp[u\ \alpha(2^{-p}d^p(x_o,x))]\,\mu(dx)<\infty
    \]
    for all $x_o\in\X$ and all $0\leq u<2.$
\end{enumerate}
\end{cor}

\proof This is Corollary \ref{res-55} with
$\chi(x)=2^{p-1}d^p(x_o,x),$ $\theta(d)=d^p$ and the following
improvement in the treatment of the inequality (\ref{eq-52}). One
can write $ \beta(d_\chi(x_o,x))\leq \gamma\circ\chi
(x_o)+\gamma\circ\chi (x)=\gamma\circ\chi (x)$ since
$\gamma\circ\chi (x_o)=0$ in this situation. As a consequence
$\max(0,2\gamma(2^{-p}t)-2\log B)$ is a transportation function,
which is a little better than its counterpart in Corollary
\ref{res-55}.
\endproof

\begin{remark}
It is known that the standard Gaussian measure $\mu$ on $\R$
satisfies $T_2$ which is the TCI with $c(x,y)=(x-y)^2$ and the
transportation function $\alpha(t)=t/2$ (see \cite{Tal96a}). As a
consequence of Corollary \ref{res-54}-b, for all $p>2,$ there is
no function $\alpha$ in $\mathcal{C}$ except $\alpha\equiv 0$
which is a transportation function for the standard Gaussian
measure and the cost function $|x-y|^p.$
\end{remark}

\subsection{Controlling the small values of $t$}

We are going to prove a general result for the behaviour of a
transportation function in the neighbourhood of zero. By a general
result, it is meant that $\mu$ is not specified. As a consequence,
it will only be shown that under the assumption that $c\leq
\chi\oplus\chi$ where $\IX e^{\delta_o \chi}\,d\mu<\infty$ for
some $\delta_o>0,$ there are tranportation functions which are
larger than some quadratic function around zero. Obtaining better
results in this direction is difficult and requires more stringent
restrictions on the reference probability measure $\mu.$

\begin{prop}\label{res-57}
Let $c$ be a cost function satisfying (\ref{eq-21}) and $c\leq
\chi\oplus\chi$ for some nonnegative measurable function $\chi$
satisfying $\IX e^{\delta_o \chi}\,d\mu<\infty$ for some
$\delta_o>0.$ Then, $\|\chi\|_\rho$ is finite and
\[
\alpha_o(t)=\left(\sqrt{t/\|\chi\|_\rho +1}-1\right)^2,\ t\geq 0
\]
is a transportation function for $c$ and $\mu.$

In particular, for all $a\geq 0$ such that $\IX
e^{a\chi}\,d\mu\leq 2,$ $t\mapsto (\sqrt{at+1}-1)^2$ is a
transportation function.
\end{prop}

Note that $(\sqrt{at+1}-1)^2=a^2t^2/4+o_{t\rightarrow
0}(t^2)=at-2\sqrt{at}+2+ o_{t\rightarrow\infty}(1).$

The Orlicz norm $\|\chi\|_\rho$ is defined at (\ref{eq-Orlicz
norm}).

\proof Because of our assumptions, we have
$\mathcal{T}_c\leq\mathcal{T}_{d_\chi},$ see (\ref{eq-dchi}).
Hence, it is enough to show that $\alpha_o$ is a transportation
function for $d_\chi.$ But this follows from  Lemma \ref{res-56}
 below and Corollary \ref{res-32}.

 The last statement follows from a simple manipulation on the
 definition of the Orlicz norm $\|\chi\|_\rho$.
\endproof

The following lemma has been used in the previous proof.
\begin{lemma}\label{res-56}
For all $\mu$ and $\nu$ in $\mathcal{P}_\chi:=\{\nu\in\PX; \IX
\chi\,d\nu<\infty\},$ we have
\[
\mathcal{T}_{d_\chi}(\mu,\nu)=\|\chi\cdot(\mu-\nu)\|_\mathrm{TV}.
\]
\end{lemma}

\proof By Kantorovich-Rubinstein's equality (\ref{eq-KR}), we have
$\mathcal{T}_{d_\chi}(\mu,\nu)=\sup\{\IX \varphi\,d (\nu-\mu);
\varphi\in B(\X), \|\varphi\|_\mathrm{Lip}\leq 1\}$ where
$\|\varphi\|_\mathrm{Lip}\leq 1$ is equivalent to
$|\varphi(x)-\varphi(y)|\leq d_\chi(x,y)$ for all $x,y.$ One can
prove without trouble (see \cite{GozPhD}) that this is equivalent
to $|\varphi(x)-a|\leq \chi(x), \forall x$ for some real $a.$
Therefore,
\begin{eqnarray*}
  \mathcal{T}_{d_\chi}(\mu,\nu)
  &=& \sup\left\{\IX \varphi\,d (\nu-\mu);
\varphi\in B(\X): |\varphi|\leq \chi\right\} \\
   &=& \sup_{k\geq 1} \sup\left\{\IX (\chi\wedge k) \theta \,d(\nu-\mu);\theta\in B(\X):  |\theta|\leq 1\right\} \\
   &=&\|\chi\cdot(\mu-\nu)\|_\mathrm{TV}
\end{eqnarray*}
which is the desired result.
\endproof

\subsection{An application: $T_1$-inequalities}

A $T_1$-inequality is a TCI with $c=d.$ Let us denote
$\mathcal{P}_d(\X)=\{\nu\in\PX; \IX d(x_*,x)\,\nu(dx)<\infty
\textrm{ for some (and therefore all) } x_*\in\X \}.$ Suppose that
$\mu$ is in $\mathcal{P}_d(\X).$  The function $\alpha$ is said to
satisfy the $T_1$-inequality for $d$ and $\mu$ if
\begin{equation}\label{eq-T1}
    \alpha(\mathcal{T}_d(\mu,\nu))\leq H(\nu\mid\mu),\ \forall \nu\in\mathcal{P}_d(\X).
\end{equation}

\begin{thm}[$T_1$-inequalities]\label{res-58}
Let $d$ be a \lsc\ metric. Suppose that $a\geq 0$ satisfies $\IX
e^{ad(x_o,x)}\,\mu(dx)\leq 2$ for some $x_o\in\X$ and that
$\gamma\in\mathcal{C}$ satisfies $\IX
e^{\gamma(d(x_1,x))}\,\mu(dx)\leq B<\infty$ for some $x_1\in\X,$
then
\begin{equation*}
    \alpha(t)=\max\Big((\sqrt{at+1}-1)^2,2\gamma(t/2)-2\log
    B\Big),\ t\geq 0
\end{equation*}
satisfies (\ref{eq-T1}).

Conversely, if a function $\alpha$ in the class $\mathcal{C}$
satisfies (\ref{eq-T1}), then
 \[
\IX \exp[u\ \alpha(d(x_*,x)/2)]\,\mu(dx)<\infty
    \]
    for all $x_*\in\X$ and all $0\leq u<2.$
\end{thm}

\proof Gathering Corollary \ref{res-54}-a, Proposition
\ref{res-57} and  the trick (\ref{eq-54}) gives us the first
statement. The second statement is a particular instance of
Corollary \ref{res-54}-b.
\endproof

Note that by Proposition \ref{res-35} we know that it is
impossible that $\alpha$ escapes from a quadratic growth at the
origin.

Theorem \ref{res-58} extends the integral criteria for the usual $T_1(C)$-inequality
in \cite{DGW03} and \cite{BV03}. Nevertheless, the control of the constant $C$ is handled more carefully
in these cited papers.

In a forthcoming paper (see the PhD manuscript \cite{GozPhD}), one of the author
has obtained the following result which is very much in the spirit of \cite{DGW03} and \cite{BV03}.

\begin{thm}Suppose that $c(x,y)=d^p(x,y)$, that $\alpha$ satisfies (\ref{eq-39}) for some $a>0$ and that $\alpha^\mc$ is unbounded on its effective domain. Then, the following statements are equivalent :
\begin{itemize}
\item There exists $b_1>0$ such that $\alpha\left(b_1\mathcal{T}_{d^p}(\nu,\mu)\right)\leq H(\nu|\mu)$ for all $\nu\in \PX$ such that $\IX d^p(x_o,x)\, \mu(dx)<\infty$
\item There exists $b_2>0$ such that $\iint_{\X^2} e^{\alpha(b_2d^p(x,y))}\,\mu(dx)\mu(dy)<+\infty$.
\end{itemize}
\end{thm}
Further details concerning the relation between $b_1$ and $b_2$ can be found in \cite{GozPhD}.

\section{Some applications: concentration of measure and deviations of empirical processes}\label{sec-applications}

In this section, we give some applications of $T_1$-inequalities.
The first application, Theorem \ref{marton} is an easy extension of a well known result of K. Marton.
The second one, Theorem \ref{dev-emp-proc} is more original and concerns the deviations of empirical processes.

In the whole section, $d$ is a metric on $\X$ which turns $(\X,d)$ into a Polish space.

\subsection{A basic lemma}

Theorem \ref{marton} and Theorem \ref{dev-emp-proc} both rely on
the following elementary lemma.

\begin{lemma}\label{basic-lemma}
Let $\mu \in \PX$ be such that $\IX d(x_o,x)\,\mu(dx)<+\infty$, for all $x_0 \in \X$, and suppose that the $T_1$- inequality
$$
\alpha\left(\mathcal{T}_d(\mu,\nu)\right)\leq H(\nu\mid\mu),\quad\forall \nu \in \PX,
$$
holds. Then, for all $1$-Lipschitz function $\varphi$, one has
\begin{equation}
\mu\left(\varphi \geq \langle \varphi,\mu\rangle + t\right)\leq e^{-\alpha(t)},\quad\forall t>0.
\end{equation}
\end{lemma}

\proof
Let $\varphi$ a $1$-Lipschitz function. For every $n\geq 1$, let us consider $\varphi_n=\varphi\vee n\wedge-n$.
According to point b. of Theorem \ref{res-BG-T1}, one has
$$
\Lambda_{\varphi_n}(s):=\log \IX e^{s\left(\varphi_n-\langle \varphi_n,\mu\rangle\right)}
\,d\mu\leq \alpha^\circledast(s),\quad\forall s\geq0.
$$
By dominating convergence, $\langle
\varphi_n,\mu\rangle\xrightarrow[n\rightarrow+\infty]{}\langle
\varphi,\mu\rangle$. Thus by Fatou's lemma, one has
$$
\Lambda_\varphi(s):=\log \IX e^{s\left(\varphi-\langle \varphi,\mu\rangle\right)}
\,d\mu\leq \alpha^\circledast(s),\quad\forall s\geq0.
$$
Now, thanks to Chebychev argument, one has for all $t\geq 0$ :
$$
\mu\left(\varphi \geq \langle \varphi,\mu \rangle +t\right)\leq \inf_{s\geq 0}
\IX e^{s(\varphi-\langle \varphi,\mu \rangle -t)}\,d\mu \leq \inf_{s\geq 0}e^{\alpha^\circledast(s)-st}=e^{-\alpha(t)}.
$$
\endproof

\subsection{$T_1$-inequalities and concentration of measure}

Let us recall that for a given probability measure $\mu$ on a
Polish space $\X$, the concentration function of $\mu$ is defined
by
$$
\theta_\mu(r)=\sup\{1-\mu(A^r) : A \text{ borel set such that } \mu(A)\geq 1/2\},\quad\forall r>0,
$$
where
\[
A^r:=\{x\in \X : d(x,A)\leq r\}.
\]
One says that $\theta$ is a concentration function for $\mu$, if there is $r_0\geq 0$ such that
$$
\theta_\mu(r)\leq \theta(r),\quad\forall r\geq r_0,
$$
or equivalently
$$
\mu(A^r)\geq 1-\theta(r),\quad \forall r\geq r_0,\quad\forall A \text{ Borel set}.
$$
Roughly speaking, the following theorem states that if $\alpha$ is
a $T_1$-transportation function for $\mu$ then $e^{-\alpha}$ is a
concentration function for $\mu$. This link between transportation
cost inequality and concentration inequality was first noticed by
K. Marton, see \cite{Mar86}. Her result extends as follows.

\begin{thm}\label{marton}
Let $\mu \in \PX$ be such that $\IX d(x_o,x)\,\mu(dx)<+\infty$ for all $x_0 \in \X$, and suppose that the $T_1$-inequality
$$
\alpha\left(\mathcal{T}_d(\mu,\nu)\right)\leq H(\nu\mid\mu),\quad\forall \nu \in \PX,
$$
holds with an unbounded $\alpha \in \mathcal{C}$. Then for all
measurable $A$ with $\mu(A)>0$, one has the following
concentration of measure inequality :
\begin{equation}\label{conc-measure}
\mu(A^r)\geq 1- e^{-\alpha(r-r_A)},\quad\forall r\geq r_A,
\end{equation}
where $r_A :=\alpha^{-1}(-\log \mu(A)).$
\end{thm}

The following proof is different from Marton's original argument. Our proof is based
on deviation arguments while Marton's one is based on transportation.
For a proof using Marton's concentration arguments see Proposition VI.81 in \cite{GozPhD}.

\proof
The function $x\mapsto d(x,A)$ is $1$-Lipschitz. Thus,
according to Lemma \ref{basic-lemma},
$$\mu(d(\cdot,A)\geq t+\langle d(\cdot,A),\mu\rangle )\leq e^{-\alpha(t)},\quad\forall t\geq 0.$$
In order to derive (\ref{conc-measure}), the only thing to do is
to show that $\langle d(\cdot,A),\mu\rangle \leq \alpha^{-1}(-\log
\mu(A))$. Let $\nu\in \PX$ be such that $\nu(A)=1$. According to
the $T_1$-inequality satisfied by $\mu$, one has
$$
\IX d(\cdot,A)\,d\mu=\IX d(\cdot,A)\,d\mu-\IX d(\cdot,A)\,d\nu\leq \mathcal{T}_d(\mu,\nu)
\leq \alpha^{-1}(H(\nu\mid\mu)).
$$
Thus,
$$
\langle d(\cdot,A),\mu \rangle\leq \alpha^{-1}\left(\inf\left\{H(\nu\mid\mu) :  \nu(A)=1\right\}\right).
$$
Let $\mu_A \in \PX$ be defined by $d\mu_A=\frac{\1_A}{\mu(A)}d\mu$ ; clearly $\mu_A(A)=1$, so

\begin{equation}\label{why-muA}
\inf\left\{H(\nu\mid\mu) :  \nu(A)=1\right\}\leq H(\mu_A\mid\mu).
\end{equation}

An easy computation yields $H(\mu_A\mid\mu)=-\log \mu(A)$.
\endproof

Note that $d(\cdot,A)$ is unbounded so that the inequality
$\IX d(\cdot,A)\,d\mu-\IX d(\cdot,A)\,d\nu\leq \mathcal{T}_d(\mu,\nu)$ needs to be justified.
Let $\pi $ be a probability on $\X^2$ with marginals $\mu$ and $\nu$, then
$\IX d(\cdot,A)\,d\mu-\IX d(\cdot,A)\,d\nu=\iint_{\X^2}d(x,A)-d(y,A)\,
\pi(dxdy)\leq \iint_{\X^2}d(x,y)\,\pi(dxdy)$. Optimizing in $\pi$ leads to the desired result.

 \medskip\textbf{Some comments.}
In Marton's approach, the probability measure $\mu_A$ plays also a
great role. Thanks to our approach, this role can be further
explained. The choice of $\mu_A$ is optimal in the sense that
(\ref{why-muA}) holds with equality:

\begin{equation}\label{muA-Iproj}
\inf\left\{H(\nu\mid\mu) :  \nu(A)=1\right\}= H(\mu_A\mid\mu).
\end{equation}

In other words, $\mu_A$ is Csisz\'ar's $I$-projection of $\mu$ on
$\{\nu\in \PX  : \nu(A)=1\}$, see \cite{Csi75,Csi84}.
\\
If $\nu$ is such that $\nu(A)=1$, one has

\begin{align*}
H(\nu\mid\mu)&=H(\nu\mid\mu_A)+\IX \log \frac{d\mu_A}{d\mu}\,d\nu\\
&=H(\nu\mid\mu_A)+\IX \log \1_A\,d\nu-\log \mu(A)\\
&=H(\nu\mid\mu_A)+H(\mu_A\mid\mu),
\end{align*}

where the last equality follows from $\IX \log \1_A\,d\nu=0$ and $H(\mu_A\mid\mu)=-\log \mu(A)$. This proves (\ref{muA-Iproj}).

\subsection{$T_1$-inequalities and deviations bounds for empirical processes.}

Lemma \ref{basic-lemma} together with the tensorization property
of Theorem \ref{res-n tenso bis} immediately implies the following

\begin{lemma}\label{basic-lemma-bis}
Let $\mu \in \PX$ be such that $\IX d(x_o,x)\,\mu(dx)<+\infty$, for all $x_0 \in \X$, and suppose that the $T_1$-inequality
$$
\alpha\left(\mathcal{T}_d(\mu,\nu)\right)\leq H(\nu\mid\mu),\quad\forall \nu \in \PX,
$$
holds. Then for all function $Z:\X^n\rightarrow \R$ which is
$1/n$-Lipschitz with respect to the metric $d^{\oplus n}$, one has
\begin{equation}
\mu^{\otimes n}\left(Z\geq \langle \mu,Z\rangle + t\right)\leq e^{-n\alpha\left(t\right)},\quad \forall t\geq0
\end{equation}
\end{lemma}

Let us consider a class $\mathcal{G}$ of $1$-Lipschitz functions
on $\X$, and $X_i$ an iid sample of law $\mu$. Let
$Z_n^{\mathcal{G}}$ be defined by

\begin{equation}\label{emprirical-process}
Z_n^{\mathcal{G}}:=\sup_{\varphi \in \mathcal{G}} \left\{\left|\frac{1}{n}
\sum_{i=1}^n\varphi(X_i)-\IX\varphi\,d\mu\right|\right\}.
\end{equation}

As $0\leq Z_n^{\mathcal{G}}=\sup_{\varphi \in
\mathcal{G}}\left\{\left|\IX \varphi\, dL_n-\IX
\varphi\,d\mu\right|\right\}\leq \mathcal{T}_d(L_n,\mu)$, one has
$Z_n^{\mathcal{G}}\in[0,+\infty[$. Further, as a supremum of
$1/n$-Lipschitz functions, the function
$$
(x_1,\ldots,x_n)\mapsto\sup_{\varphi \in \mathcal{G}}
\left\{\left|\frac{1}{n}\sum_{i=1}^n\varphi(x_i)-\IX\varphi\,d\mu\right|\right\}
$$
is $1/n$-Lipschitz too. This implies in particular that
$Z_n^{\mathcal{G}}$ is measurable. The random variable
$Z_n^{\mathcal{G}}$ is called an empirical process. Applying Lemma
\ref{basic-lemma-bis}, one immediately obtains the following
theorem.

\begin{thm}\label{dev-emp-proc}
Let $\mu \in \PX$ be such that $\IX d(x_o,x)\,\mu(dx)<+\infty$,
for all $x_0 \in \X$, and suppose that the $T_1$-inequality
$$
\alpha\left(\mathcal{T}_d(\mu,\nu)\right)\leq H(\nu\mid\mu),\quad\forall \nu \in \PX,
$$
holds. If $\mathcal{G}$ is a class of $1$-Lipschitz functions on
$\X$ then the empirical process $Z_n^{\mathcal{G}}$ defined by
(\ref{emprirical-process}) satisfies the following inequality
\begin{equation}
\P\left(Z_n^{\mathcal{G}}\geq\E\left[Z_n^{\mathcal{G}}\right]+t\right)\leq e^{-n\alpha(t)},\quad\forall t\geq0.
\end{equation}
\end{thm}

The literature about the deviations of empirical processes is
huge. For a good overview of this subject, one can read
P.~Massart's Saint-Flour lecture notes \cite{Mass}.

Now, if $(\X,\|\,\cdot\,\|)$ is a Banach space, and $\mu \in \PX$
such that $\IX \|x\|\,d\mu<+\infty$ then taking
$\mathcal{G}=\{\ell \in \X^* : \|\ell\|_{\X^*}=1 \}$, where $\X^*$
is the topological dual space of $\X$,  one obtains
$$
Z_n^{\mathcal{G}}=\left\|\frac{1}{n}\sum_{i=1}^nX_i-\IX x\,d\mu\right\|,
$$
where $\IX x\,\mu(dx)$ is well defined in the Bochner sense. In this special case, we have the following result.

\begin{thm}\label{dev-moy-emp}
Let $\mu \in \PX$ be such that $\IX \|x\|\,\mu(dx)<+\infty$, and suppose that the $T_1$-inequality
$$
\alpha\left(\mathcal{T}_{\|\,\cdot\,\|}(\mu,\nu)\right)\leq H(\nu\mid\mu),\forall \nu \in \PX,
$$
holds. If $X_i$ is an iid sequence of law $\mu$, then letting
$Z_n=\left\|\frac{1}{n}\sum_{i=1}^nX_i-\IX x\,d\mu\right\|$, one
has
\begin{equation}
\P\left(Z_n\geq\E\left[Z_n\right]+t\right)\leq e^{-n\alpha(t)},\quad\forall t\geq0.
\end{equation}
\end{thm}

\begin{remark} In order to obtain precise deviations results for
$Z_n^{\mathcal{G}}$ (resp. $Z_n$), one must be able to estimate
the term $\E\left[Z_n^{\mathcal{G}}\right]$ (resp.
$\E\left[Z_n\right]$).
\end{remark}

Let us give some examples.
\par\medskip
\textbf{Example 1. Quantitative versions of Sanov theorem.}
Suppose that $\mathcal{G}$ is the set of all bounded $1$-Lipschitz
functions on $\X$, then
$Z_n^{\mathcal{G}}=\mathcal{T}_d(L_n,\mu),$ see (\ref{eq-KR}).
\\
The following theorem is Theorem 10.2.1 of \cite{RacRus} (volume II).

\begin{thm}\label{Rac-Rusch}
Let $\mu$ be a probability measure on $\R^q$ (equipped with its usual euclidean norm $\|\,\cdot\,\|_2$) such that
\begin{equation}\label{Rac-Rusch1}
c:=\int \|x\|_2^{q+5}\,d\mu<+\infty.
\end{equation}
Then, there is $D>0$ depending only on $c$ and $q$, such that
\begin{equation}\label{Rac-Rusch2}
\E\left[\mathcal{T}_{d_2}(L_n,\mu)\right]\leq Dn^{-\frac{1}{q+4}},
\end{equation}
where $d_2$ is the metric associated to $\|\,\cdot\,\|_2$.
\end{thm}

Thanks to this result, one obtains the following quantitative
version of Sanov theorem :
\begin{cor}
Let $\mu$ be a probability on $\R^q$, satisfying (\ref{Rac-Rusch1}) and the $T_1$-inequality
$$
\alpha\left(\mathcal{T}_{d_2}(\mu,\nu)\right)\leq H(\nu\mid\mu),\quad\forall \nu \in \mathcal{P}(\R^q),
$$
where $d_2$ is the usual euclidean metric on $\R^q$. Then, the following inequality holds :
$$
\P\left(\mathcal{T}_{d_2}(L_n,\mu)\geq t\right)
\leq \exp\left(-n\alpha\left(t-\frac{D}{n^{\frac{1}{q+4}}}\right)\right),\quad\forall t>0,\quad\forall n\geq \left(\frac{D}{t}\right)^{q+4},
$$
where $D$ is the constant of (\ref{Rac-Rusch2}).
\end{cor}

In \cite{BGV05}, F. Bolley, A. Guillin and C. Villani have also
obtained a quantitative version of Sanov theorem with alternative
arguments.

\par\medskip
\textbf{Example 2. Deviations bounds for empirical means.}
Let $\X$ be a Banach space and consider
\begin{equation}\label{Z_n}
Z_n=\left\|\frac{1}{n}\sum_{i=1}^nX_i-\IX x\,d\mu\right\|,
\end{equation}
where $X_i$ is an iid sequence of law $\mu$. In order to control
the term $\E[Z_n]$, a classical assumption is to require that $\X$
is of type $p>1$, \textit{ie} there is $b>0$ such that for every
sequence $(Y_i)_i$ of centered random variables with
$\E\left[\|Y_i\|^p\right]<+\infty$, one has

\begin{equation}\label{type-p}
\E\left[\left\|Y_1+\cdots +Y_n\right\|^p\right]\leq b \left[\E\left[\left\|Y_1\right\|^p\right]
+\cdots +\E\left[\left\|Y_n\right\|^p\right]\right].
\end{equation}

If $\X$ is of type $p$ and $\E\left[\|X_1\|^p\right]<+\infty$,
then one can deduce immediately from (\ref{type-p}) the following
control:

\begin{equation}\label{type-p2}
\E\left[Z_n\right]\leq \frac{1}{n^{1-1/p}}\left(b\E\left[\|X_1-\E[X_1]\|^p\right]\right)^{1/p}.
\end{equation}

Controls like (\ref{type-p2}) can be used in Theorem
\ref{dev-moy-emp} to derive precise deviations bounds for
empirical means. Let us conclude this section with a concrete
example.

\begin{thm}
Let $\mu$ be a probability measure on a Banach space
$(\X,\|\,\cdot\,\|)$ such that $\IX e^{a\|x\|}\,\mu(dx)<+\infty$,
for some $\delta>0$. Then, for all sequence $X_i$ of iid random
variables with law $\mu$, one has
\begin{equation}\label{Yur-like}
\P\left(Z_n\geq \E[Z_n]+t\right)\leq e^{-n\left(\sqrt{1+\frac{t}{M}}-1\right)^2},\quad\forall t>0,
\end{equation}
where $Z_n$ is defined by (\ref{Z_n}) and $M:=\inf\left\{b>0 : \iint_{\X^2}e^{\frac{\|x-y\|}{b}}\mu(dx)\mu(dy)\leq2\right\}.$
\end{thm}

\proof According to Corollary \ref{res-33}, $\mu$ satisfy the
$T_1$-inequality
$$
\alpha\left(\mathcal{T}_{\|\,\cdot\,\|}(\mu,\nu)\right)\leq H(\nu\mid\mu),\quad\forall \nu \in \PX,
$$
with $\alpha(t)=\left(\sqrt{1+\frac{t}{M}}-1\right)^2$. Thus,
applying Theorem \ref{dev-moy-emp}, the result follows
immediately.
\endproof

Inequality (\ref{Yur-like}) is very close to a well known
inequality by Yurinskii (\cite{Yur76}, Theorem 2.1). Under the
same assumptions on $\mu$, one can easily derive from Yurinskii's
result the following bound :
\begin{equation}\label{Yur}
\P\left(Z_n\geq \E\left[Z_n\right]+t\right)\leq \exp\left(-\frac{1}{8}\frac{nt^2}{2M_0^2+tM_0}\right),\quad\forall t >0,
\end{equation}
where $M_0=\inf\left\{b>0 : \IX e^{\frac{\|x\|}{b}}\,\mu(dx)\leq 2\right\}$. To compare (\ref{Yur-like}) and (\ref{Yur}) first note that
\begin{equation}\label{comparaison1}
\left(\sqrt{1+u}-1\right)^2\geq \frac{u^2}{2(2+u)},\quad\forall u>0,
\end{equation}
(this is left to the reader). Next, let us show that
\begin{equation}\label{comparaison2}
M\leq 2M_0.
\end{equation}
This follows from the following inequality :
$$
\iint_{\X^2} e^{\frac{\|x-y\|}{2M_0}}\,\mu(dx)\mu(dy)\stackrel{(i)}{\leq} \left(\IX e^{\frac{\|x\|}{2M_0}}\,\mu(dx)\right)^2
\stackrel{(ii)}{\leq}\IX e^{\frac{\|x\|}{M_0}}\,\mu(dx)\stackrel{(iii)}{\leq} 2,$$
where (i) comes from the triangle inequality, (ii) from Jensen inequality and (iii) from the definition of $M_0$.
Thanks to (\ref{comparaison1}) and (\ref{comparaison2}), one obtains
$$\left(\sqrt{1+\frac{t}{M}}-1\right)^2\geq\frac{t^2}{2(2M^2+tM)}\geq\frac{t^2}{8(2M_0^2+tM_0/2)}\geq\frac{t^2}{8(2M_0^2+tM_0)}.$$
Thus, (\ref{Yur-like}) is a little bit stronger than (\ref{Yur}).\\
Yurinskii's proof relies on martingale arguments, while our proof is a direct consequence of the tensorization mechanism.

\section{Large deviations and $\mathcal{T}$-inequalities. Abstract
results}\label{sec-abstract results}

The framework is the same as in Section \ref{sec-integral criteria}. See in particular Remark \ref{rem-51}.

\subsection{A deviation function is a transportation function}

In this section, we give a rigorous proof at Theorem
\ref{res-abstract} of the Recipe \ref{recipe2}  for an increasing
deviation function which may possibly be not convex. This extends
Proposition \ref{res-51}.

\begin{thm}\label{res-abstract}
Let us assume (\ref{hyp-F}) and (\ref{hyp-Phi}).
\begin{enumerate}
    \item[(a)]
    Any deviation function is a
transportation function.
    \item[(b)]
If in addition $\mathcal{T}$ is continuous on $\NF,$ then the
converse also holds: any transportation function is a deviation
function.
\end{enumerate}
\end{thm}

\proof (a)\quad As $\mathcal{T}$ is \lsc, for all $t\geq 0$ the
set $\{\nu\in\NF; \mathcal{T}(\nu)>t\}$ is open. It follows with
the LD lower bound that
\begin{equation*}
  -\inf\{H(\nu\mid\mu); \nu\in\NF, \mathcal{T}(\nu)>t\}
   \leq \loginf \P(\mathcal{T}(L_n)>t)
\end{equation*}
Let $\alpha$ be any deviation function: for all $t\geq 0,$
$\logsup\P(\mathcal{T}(L_n)\geq t)\leq -\alpha(t).$ Hence we
obtain $\alpha(t)\leq \inf\{H(\nu\mid\mu); \nu\in\NF,
\mathcal{T}(\nu)>t\}$ so that $\alpha(t-\delta)\leq H(\nu\mid\mu)$
for all $\nu\in\NF$ and $\delta>0$ such that $\mathcal{T}(\nu)>
t-\delta.$ Taking $t=\mathcal{T}(\nu)$ leads us to
$\alpha(\mathcal{T}(\nu)-\delta)\leq H(\nu\mid\mu)$  for all
$\nu\in\NF$ and $\delta>0.$ As $\alpha$ is increasing and
$\delta>0$ is arbitrary, we have $\alpha(\mathcal{T}(\nu)^-)\leq
H(\nu\mid\mu).$ The desired result follows from the assumed left
continuity of $\alpha.$
\par\medskip\noindent
(b)\quad As $\mathcal{T}$ is continuous, because of the
contraction principle, $\{\mathcal{T}(L_n)\}$ obeys the LDP with
rate function $i(t)=\inf\{H(\nu\mid\mu); \nu\in\NF,
\mathcal{T}(\nu)=t\}, t\geq 0.$ In particular, the LD upper bound:
$\logsup \P(\mathcal{T}(L_n)\geq t)\leq -\inf\{i(s); s\geq t\},$
is satisfied.
\\
Let $\alpha$ be a transportation function. It clearly satisfies
$\alpha(t)\leq \inf\{H(\nu\mid\mu); \nu\in\NF,
\mathcal{T}(\nu)=t\}$ for all $t.$ That is: $\alpha\leq i.$
Finally, for all $t\geq 0,$
\begin{eqnarray*}
  \logsup \P(\mathcal{T}(L_n)\geq t) &\leq & -\inf_{s\geq t}i(s) \\
   &\leq& -\inf_{s\geq t}\alpha (s) \\
   &=& -\alpha(t)
\end{eqnarray*}
where the last equality holds because $\alpha$ is increasing. This
means that $\alpha$ is a deviation function.
\endproof

\textit{Remarks}.
\begin{itemize}
    \item Note that we didn't use the specific form (\ref{def-T}) of
$\mathcal{T},$ but only its lower semicontinuity.
    \item Similarly, we didn't use the specific properties of the
    relative entropy, but only that it is a LDP rate function for
    $\{L_n\}.$
    \item Statement (b) will not be used later, but it is satisfactory to
    know that a transportation function is not far from being a deviation
    function. A natural situation where $\mathcal{T}$ is
    continuous appears with $c=d^p$ since the Wasserstein's metric
    $\mathcal{T}_{d^p}^{1/p}$ metrizes $\sigma(\NF,\mathcal{F})$
    with $\mathcal{F}$ the space of all continuous functions
    $\varphi$ such that $|\varphi(x)|\leq c(1+d(x_o,x)^p), \forall
    x$ for some constant $c,$ see (\cite{Vill}, Chapter 7).
\end{itemize}

\subsection{The transportation function $J_\Phi$}
\label{sec-JPhi}

With Theorem \ref{res-abstract} in hand,
it is enough to compute a deviation function $\alpha$ to obtain
the TCI
\begin{equation}\label{eq-TIneq2}
\alpha(\mathcal{T}(\nu)) \leq H(\nu\mid\mu),\quad \forall\nu\in\NF
\end{equation}

But these functions may be rather hard to compute because of the
sup in the definition (\ref{def-T}) of
$$
T_n=\mathcal{T}(L_n)=\sup_{(\psi,\varphi)\in\Phi}\{\langle\varphi,L_n\rangle+\langle\psi,\mu\rangle\}.
$$
However, it is shown at Theorem \ref{res-reassuring} below, that
more can be said about transportation functions.
\\

\par\medskip\noindent \textbf{Assumptions (A).}\
\begin{it}
The following requirements are assumed to hold.

\begin{enumerate}
    \item[(i)]
    We assume (\ref{hyp-F}):
    \begin{equation*}
    \IX e^\varphi\,d\mu<\infty, \forall
    \varphi\in \mathcal{F}.
    \end{equation*}
    \item[(ii)]
        We assume (\ref{hyp-Phi}):
    \begin{equation*}
   (0,0)\in \Phi\subset \mathcal{F}\times\mathcal{F},
    \end{equation*}
     \item[(iii)]
    For all $(\psi,\varphi)\in\Phi,$ $\psi+\varphi\leq 0.$\\
\end{enumerate}
\end{it}

Requirement (iii) always holds in the norm case: $\Phi=\Phi_U,$
and it holds in the transportation case $\Phi=\Phi_c$ if
$c(x,x)=0, \forall x\in\X$.

Let us define
    \begin{equation*}
    \Lambda(\varphi):=\log \IX e^\varphi\,d\mu.
    \end{equation*}

\begin{prop}\label{res-LDP-a}
Under the assumption (\ref{hyp-F})
\begin{enumerate}
    \item[(a)]
    $\{L_n\}$ obeys the LDP in $\NF$ with the rate function
    \begin{equation}\label{Lambda star}
    H(\nu\mid\mu)= \Lambda^*(\nu)=\sup_{\varphi\in\mathcal{F}}\{\langle\varphi,\nu\rangle-\Lambda(\varphi)\}, \quad
    \nu\in\NF.
    \end{equation}
    \item[(b)]
    and for all $(\psi,\varphi)\in\Phi,$ $\{T_n^{\psi,\varphi}\}_{n\geq 1}$ obeys the LDP in $\mathbb{R}$
    with the rate function
    \begin{equation*}
    J_{\psi,\varphi}(t)=\sup_{s\in \mathbb{R}}\{st-\Lambda(s\varphi)-s\langle\psi,\mu\rangle\},
    \ t\in\mathbb{R}.
    \end{equation*}

\end{enumerate}
\end{prop}
\proof Statement (a) is Theorem \ref{res-Sanov}.
\\
The function $J_{\psi,\varphi}$ is the convex conjugate of
    $$\Lambda_{\psi,\varphi}(s):=\Lambda(s\varphi)+s\langle\psi,\mu\rangle,\
    s\in\mathbb{R}.$$
Since $\Lambda_{\psi,\varphi}$ is a steep function under
assumptions (ii) and (iii), (b) is a direct consequence of
G\"artner-Ellis theorem.
\endproof

We know that $J_{\psi,\varphi}$ is convex with a minimum value 0
attained at $\Lambda_{\psi,\varphi}^{\prime}(0).$ Under assumption
(iii), we have
$\Lambda_{\psi,\varphi}^{\prime}(0)=\langle\varphi+\psi,\mu\rangle\leq
0.$ Therefore, $J_{\psi,\varphi}$ is an increasing nonnegative
function on $[0,\infty)$ and so are  $J_\Phi$ and
$\widetilde{J}_\Phi$ given by
\begin{eqnarray}\label{eq-JPhi}
  J_\Phi(t)&:=&\widetilde{J}_\Phi(t^-), t>0 \mathrm{\quad where} \\
  \widetilde{J}_\Phi(t)&:=&\inf_{(\psi,\varphi)\in\Phi}J_{\psi,\varphi}(t)\in[0,\infty],
t\geq 0\notag
\end{eqnarray}
with $J_\Phi(0)=0.$ This last equality follows from assumption
(ii). As $\Lambda_{\psi,\varphi}^{\prime}(0)\leq 0,$ it also holds
that for all $t\geq 0,$
$J_{\psi,\varphi}(t)=\Lambda^\mc_{\psi,\varphi}(t):=\sup_{s\geq
0}\{st-\Lambda_{\psi,\varphi}(s)\}$ where the sup is taken over
$s\geq 0$ rather than $s\in\mathbb{R}.$ It follows that one can
equivalently define $J_\Phi$ as follows.

\begin{definition}[of the functions $J_\Phi$ and $J$].\label{def-J}
\begin{itemize}
    \item $J_\Phi$ is the left continuous version of the
    increasing function
\begin{equation*}
    t\in[0,\infty)\mapsto\inf_{(\psi,\varphi)\in\Phi}\sup_{s\geq 0}
    \{st-\Lambda(s\varphi)-s\langle\psi,\mu\rangle\}\in [0,\infty].
\end{equation*}
    \item $J$ is the \emph{best transportation function}. Clearly, it is
the left continuous function of the increasing function
\begin{equation*}
    t\in[0,\infty)\mapsto\inf\{H(\nu\mid\mu); \nu\in\NF: \mathcal{T}(\nu)\geq t\}\in [0,\infty].
\end{equation*}
\end{itemize}
\end{definition}

Although the best transportation function $J$ might be out of
reach in many situations, we have the following reassuring result.

\begin{thm}\label{res-reassuring}
Suppose that Assumptions (A) hold. Then, $J_\Phi$ is a
transportation function and the \emph{best} transportation
function \emph{in the class $\mathcal{C}$} is the convex \lsc\
regularization of $J_\Phi.$
\end{thm}

\proof This statement is a collection of the statements of Theorem
\ref{res-JPhi}-a and Corollary \ref{res-C-transport}-a,b which
will be proved below.
\endproof

\begin{thm}\label{res-JPhi}
Suppose that Assumptions (A) hold.
\begin{enumerate}
    \item[(a)]
    Then, $J_\Phi$ is a transportation function for $\mathcal{T}$
    and
    $\{L_n\}.$ This can be equivalently rewritten as the following TCI
     \begin{equation*}
     J_\Phi(\mathcal{T}(\nu))\leq H(\nu\mid\mu),\
     \forall\nu\in\NF.
     \end{equation*}
    \item[(b)]
    If in addition $\mathcal{T}$ is continuous on $\NF,$ then
    $J_\Phi$ is the best transportation function. It
    is also the best deviation function: This means that
    $J_\Phi=J.$
\end{enumerate}
\end{thm}
\proof (a)\quad As
$\nu\mapsto\langle\varphi,\nu\rangle+\langle\psi,\mu\rangle$ is
continuous, it follows from the contraction principle that
$J_{\psi,\varphi}(t)=\inf\{H(\nu\mid\mu); \nu\in\NF,
\langle\varphi,\nu\rangle+\langle\psi,\mu\rangle=t\}$ for all
$t\geq 0.$ Hence,
$J_{\psi,\varphi}(\langle\varphi,\nu\rangle+\langle\psi,\mu\rangle)\leq
H(\nu\mid\mu)$ for all $\nu\in\NF$ and a fortiori
\begin{equation*}
    \widetilde{J}_\Phi(\langle\varphi,\nu\rangle+\langle\psi,\mu\rangle)\leq
    H(\nu\mid\mu),
\end{equation*}
as soon as $\langle\varphi,\nu\rangle+\langle\psi,\mu\rangle\geq
0.$ As $\widetilde{J}_\Phi$ is increasing, by the definition
(\ref{def-T}) of $\mathcal{T}(\nu),$ one obtains:
$\widetilde{J}_\Phi(\mathcal{T}(\nu)^-)\leq H(\nu\mid\mu)$ which
is the desired result. Note that $\mathcal{T}(\nu)\geq 0$ since
$(0,0)\in \Phi$ (assumption (A.ii)).
\par\medskip\noindent
(b)\quad Because of part (b) of Theorem \ref{res-abstract}, it is
enough to prove that $J_\Phi=J.$ Because of part (a) of the
present theorem, $J_\Phi$ is a transportation function, and by
part (b) of Theorem \ref{res-abstract}, it is also a deviation
function. Therefore, $J_\Phi\leq J$ and it remains to prove that
$J\leq J_\Phi.$
\\
By the LD lower bound for $\{T^{\psi,\varphi}_n\},$ for all $t\geq
0,$
\begin{eqnarray*}
  -\inf_{r>t}J_{\psi,\varphi}(r) &\leq& \loginf \P(T^{\psi,\varphi}_n>t) \\
   &\leq& \logsup \P\left(\sup_{(\psi,\varphi)\in\Phi}T^{\psi,\varphi}_n\geq t\right) \\
   &\leq& -J(t).
\end{eqnarray*}
Since $J_{\psi,\varphi}$ is increasing, we have: $J(t)\leq
\inf_{r>t}J_{\psi,\varphi}(r)=J_{\psi,\varphi}(t^+),$ so that for all $t\geq 0$
\begin{align*}
    J(t)&\leq \inf\{J_{\psi,\varphi}(t^+), (\psi,\varphi)\in\Phi\}\\
    &= \inf_\phi \inf_{u>t} J_\phi(u)\\
    &= \inf_{u>t} \inf_\phi J_\phi(u)\\
    &= \widetilde{J}_\Phi (t^+).
\end{align*}

 As $J$ and
$\widetilde{J}_\Phi$ are increasing and $J$ is left continuous,
this gives $J(t)\leq \widetilde{J}_\Phi(t^-)$ for all $t>0$ which
is the desired result.
\endproof

\subsection{Connections with Theorem \ref{res-BG}}
Let us first give an alternative proof of
 criterion $(b)\Rightarrow (a)$ of Theorem \ref{res-BG}.
\\
 We keep the Assumptions (A) of
Section \ref{sec-JPhi}. Note that because of Assumptions (A.ii)
and (A.iii), the function
\begin{equation}\label{def-LambdaPhi}
    \Lambda_\Phi(s):=\sup_{(\psi,\varphi)\in\Phi}\Lambda_{\psi,\varphi}(s)
    =\sup_{(\psi,\varphi)\in\Phi}\{\Lambda(s\varphi)+s\langle\psi,\mu\rangle\},
    \ s\geq 0
\end{equation}
 is in the class $\mathcal{C}.$
It follows that its monotone conjugate
$$
\Lambda_\Phi^\mc(t)=\sup_{s\geq 0}\{st-\Lambda_\Phi(s)\}, t\geq 0
$$
is also in $\mathcal{C}.$  Thanks to formula (\ref{eq-JPhi}), for
all $t\geq 0,$ we have
\begin{eqnarray*}
  \Lambda_\Phi^\mc(t)&\leq& \sup_{s\geq
    0}\left\{st-\sup_{(\psi,\varphi)\in\Phi}\Lambda_{\psi,\varphi}(s)\right\} \\
    &=& \sup_{s\geq
    0}\inf_{(\psi,\varphi)\in\Phi}\{st-\Lambda_{\psi,\varphi}(s)\}\\
   &\leq& \inf_{(\psi,\varphi)\in\Phi}\sup_{s\geq
   0}\{st-\Lambda_{\psi,\varphi}(s)\}\\
   &=& \widetilde{J}_\Phi(t)
\end{eqnarray*}
But $\Lambda_\Phi^\mc(t)$ is left continuous, hence
\begin{equation}\label{eq:b}
    \Lambda_\Phi^\mc\leq J_\Phi.
\end{equation}
As $J_\Phi$ is a transportation function (Theorem \ref{res-JPhi}),
so is $\Lambda_\Phi^\mc.$

The criterion $(b)\Rightarrow (a)$ of Theorem \ref{res-BG} follows
from the above considerations. Indeed, (b) states that
$\Lambda_\Phi\leq \alpha^\mc.$ Therefore, with (\ref{eq:b}):
$\alpha\leq \Lambda_\Phi^\mc\leq J_\Phi.$ Hence, $\alpha$ is a
transportation function.

An easy consequence of Theorem \ref{res-BG} is the following

\begin{cor}\label{res-C-transport}
Suppose that Assumptions (A) hold.
\begin{itemize}
    \item[(a)]
    The best transportation function in the class $\mathcal{C}$ is
$\Lambda_\Phi^\mc.$ This means that $\alpha\in\mathcal{C}$ is a
transportation function if and only if $\alpha\leq
\Lambda_\Phi^\mc.$
    \item[(b)]
    Moreover, $\Lambda^\mc_\Phi$ is the convex \lsc\ regularization of
    $J_\Phi$ (in restriction to $t\in [0,\infty)$).
    \item[(c)]
    If $\mathcal{T}$ is continuous, then $\Lambda^\mc_\Phi$ is also
    the best deviation function in the class $\mathcal{C}.$
\end{itemize}
\end{cor}

\proof
 The best function $\alpha^\mc\in\mathcal{C}$ satisfying (b)
of Theorem \ref{res-BG} is $\alpha^\mc=\Lambda_\Phi,$ see
(\ref{def-LambdaPhi}). Because of the equivalence
$(a)\Leftrightarrow (b)$ of Theorem \ref{res-BG}, its monotone
conjugate $\Lambda^\mc_\Phi$ is the best transportation function
in $\mathcal{C}.$ This is (a).

Let us prove (b). In order to work with usual convex conjugates,
let us state $J_\phi(t)=+\infty$ for all $t<0$ and $\phi\in\Phi.$
We have
\begin{eqnarray*}
  (\inf_\phi J_\phi)^*(s) &=& \sup_t \{st-\inf_\phi J_\phi(t)\} \\
   &=& \sup_{t,\phi}\{st-J_\phi(t)\}\\
   &=& \sup_{\phi}\sup_{t}\{st-J_\phi(t)\} \\
   &=& \sup_{\phi} J_\phi^*(s).
\end{eqnarray*}
Hence, the convex \lsc\ regularization  of $J_\Phi:=\inf_\phi
J_\phi$ is $(\inf_\phi J_\phi)^{**}=(\sup_{\phi}
J_\phi^*)^*=(\sup_{\phi} \Lambda_\phi^{**})^*$ But, the convex
\lsc\ regularization of $\sup_\phi \Lambda_\phi$ is $\sup_{\phi}
\Lambda_\phi^{**}.$ Therefore, $J_\Phi^{**}=(\sup_{\phi}
\Lambda_\phi^{**})^*=(\sup_{\phi} \Lambda_\phi)^*=\Lambda_\Phi^*.$
But it is already seen that in restriction to $t\in[0,\infty),$
$\Lambda_\Phi$ is in $\mathcal{C},$ so that
$\Lambda_\Phi^*(t)=\Lambda_\Phi^\mc(t)$ for all $t\geq 0.$
\\
Finally, (c) is a direct consequence of (b) and Theorem
\ref{res-JPhi}-(b).
\endproof


\begin{thebibliography}{10}

\bibitem{BGL01}
S.~G. Bobkov, I.~Gentil, and M.~Ledoux.
\newblock Hypercontractivity of {H}amilton-{J}acobi equations.
\newblock {\em Journal de {M}ath\'ematiques {P}ures et {A}plliqu\'ees},
  80(7):669--696, 2001.

\bibitem{BG99}
S.G. Bobkov and F.~G\"{o}tze.
\newblock Exponential integrability and transportation cost related to
  logarithmic {S}obolev inequalities.
\newblock {\em Journal of Functional Analysis.}, 163:1--28, 1999.

\bibitem{BGV05}
F.~Bolley, A.~Guillin, and C.~Villani.
\newblock Quantitative concentration inequalities for empirical measures on
  non-compact spaces.
\newblock preprint.\\ Available online via {\tt
  http://www.ceremade.dauphine.fr/\~{}guillin/index3.html}, 2005.

\bibitem{BV03}
F.~Bolley and C.~Villani.
\newblock Weighted {C}sisz\'ar-{K}ullback-{P}insker inequalities and
  applications to transportation inequalities.
\newblock To appear in Annales de la Facult\'e des Sciences de Toulouse.
  Available online via {\tt
  http://www.umpa.ens-lyon.fr/\~{}cvillani/cv.html\#publicationlist}, 2005.

\bibitem{Csi75}
I.~Csisz\'ar.
\newblock {$I$}-divergence geometry of probability distributions and
  minimization problems.
\newblock {\em Annals of Probability}, 3:146--158, 1975.

\bibitem{Csi84}
I.~Csisz\'ar.
\newblock Sanov property, generalized {$I$}-projection and a conditional limit
  theorem.
\newblock {\em Annals of Probability}, 12:768--793, 1984.

\bibitem{DZ}
A.~Dembo and O.~Zeitouni.
\newblock {\em Large deviations techniques and applications. {S}econd edition}.
\newblock Applications of Mathematics 38. Springer Verlag, 1998.

\bibitem{DGW03}
H.~Djellout, A.~Guillin, and L.~Wu.
\newblock Transportation cost-information inequalities for random dynamical
  systems and diffusions.
\newblock {\em Annals of Probability}, 32(3B):2702--2732, 2004.

\bibitem{EiS}
P.~Eichelsbacher and U.~Schmock.
\newblock Large deviations of {U}-empirical measures in strong topologies and
  applications.
\newblock {\em Annales de l'Institut Henri Poincar\'e}, 38(5):779--797, 2002.

\bibitem{GozPhD}
N.~Gozlan.
\newblock Principe conditionnel de {G}ibbs pour des contraintes fines
  approch\'ees et in\'egalit\'es de transport.
\newblock PhD Thesis, Universit\'e de Paris 10, 2005.

\bibitem{Led01}
M.~Ledoux.
\newblock {\em The Concentration of Measure Phenomenon}.
\newblock Mathematical Surveys and Monographs 89. American Mathematical
  Society, Providence RI, 2001.

\bibitem{LeoN02}
C.~L\'eonard and J.~Najim.
\newblock An extension of {S}anov's theorem : application to the {G}ibbs
  conditioning principle.
\newblock {\em Bernoulli}, 8(6):721--743, 2002.

\bibitem{Mar86}
K.~Marton.
\newblock A simple proof of the blowing-up lemma.
\newblock {\em IEEE Transactions on Information Theory}, 32:445--446, 1986.

\bibitem{Mar96}
K.~Marton.
\newblock Bounding $\bar{d}$-distance by informational divergence: a way to
  prove measure concentration.
\newblock {\em Annals of Probability}, 24:857--866, 1996.

\bibitem{Mass}
P.~Massart.
\newblock Saint-{F}lour {L}ecture {N}otes.
\newblock \\ Available online via {\tt http://www.math.u-psud.fr/\~{}massart/},
  2003.

\bibitem{OVill00}
F.~Otto and C.~Villani.
\newblock Generalization of an inequality by {T}alagrand and links with the
  logarithmic {S}obolev inequality.
\newblock {\em Journal of Functional Analysis}, 173:361--400, 2000.

\bibitem{RacRus}
S.~Rachev and L.~R\"uschendorf.
\newblock {\em Mass Transportation Problems. Vol I : Theory, Vol. II :
  Applications}.
\newblock Probability and its applications. Springer Verlag, New York, 1998.

\bibitem{Tal96a}
M.~Talagrand.
\newblock Transportation cost for gaussian and other product measures.
\newblock {\em Geometric and Functional Analysis}, 6:587--600, 1996.

\bibitem{Vill}
C.~Villani.
\newblock {\em Topics in Optimal Transportation}.
\newblock Graduate Studies in Mathematics 58. American Mathematical Society,
  Providence RI, 2003.

\bibitem{Yur76}
V.V. Yurinskii.
\newblock Exponential inequalities for sums of random vectors.
\newblock {\em Journal of multivariate analysis}, 6:473--499, 1976.

\end{thebibliography}
\end{document}